\documentclass[a4paper,11pt,reqno,oneside]{amsart}
\usepackage{amsmath,amssymb,amsthm,mathtools}
\usepackage{enumitem,xcolor}
\usepackage[scale=0.75,hcentering]{geometry}
\usepackage{cite,xurl,microtype,aliascnt,hyperref}
\usepackage[capitalize,nameinlink,noabbrev]{cleveref}
\makeatletter
\@ifundefined{subjclassname@2020}{\@namedef{subjclassname@2020}{%
\textup{2020} Mathematics Subject Classification}}{}
\makeatother
\definecolor{cobalt}{rgb}{0.0,0.28,0.67}
\hypersetup{bookmarksopen=true,
pdftitle={Rank Additivity and Comon's Conjecture for Real Hankel Tensors},
pdfauthor={Jianze Li},
pdfsubject={Rank Additivity and Comon's Conjecture for Real Hankel Tensors},
pdfkeywords={Hankel tensor, real rank, symmetric rank, partially symmetric rank, direct sum, apolarity},
colorlinks=true,linkcolor=cobalt,citecolor=cobalt,urlcolor=cobalt}
\theoremstyle{plain}
\newtheorem{theorem}{Theorem}[section]
\newaliascnt{proposition}{theorem}
\newtheorem{proposition}[proposition]{Proposition}
\aliascntresetthe{proposition}
\newaliascnt{lemma}{theorem}
\newtheorem{lemma}[lemma]{Lemma}
\aliascntresetthe{lemma}
\newaliascnt{corollary}{theorem}
\newtheorem{corollary}[corollary]{Corollary}
\aliascntresetthe{corollary}
\newaliascnt{problem}{theorem}
\newtheorem{problem}[problem]{Open Problem}
\aliascntresetthe{problem}
\newaliascnt{proofobligation}{theorem}

\aliascntresetthe{proofobligation}
\theoremstyle{remark}
\newaliascnt{remark}{theorem}
\newtheorem{remark}[remark]{Remark}
\aliascntresetthe{remark}
\crefname{theorem}{Theorem}{Theorems}
\crefname{proposition}{Proposition}{Propositions}
\crefname{lemma}{Lemma}{Lemmas}
\crefname{corollary}{Corollary}{Corollaries}
\crefname{remark}{Remark}{Remarks}
\crefname{problem}{Open Problem}{Open Problems}
\crefname{proofobligation}{Proof obligation}{Proof obligations}
\newcommand{\tens}[1]{\boldsymbol{#1}}
\newcommand{\matr}[1]{\boldsymbol{#1}}
\newcommand{\vect}[1]{\boldsymbol{#1}}
\newcommand{\set}[1]{\mathcal{#1}}
\newcommand{\eqdef}{\mathrel{\triangleq}}
\newcommand{\RR}{\mathbb R}
\newcommand{\CC}{\mathbb C}
\newcommand{\rk}{r}
\DeclareMathOperator{\Sym}{\mathsf{Sym}}
\DeclareMathOperator{\Span}{\mathsf{span}}
\DeclareMathOperator{\codim}{\mathsf{codim}}

\DeclareMathOperator{\supp}{\mathsf{supp}}
\DeclareMathOperator{\lcm}{\mathsf{lcm}}
\numberwithin{equation}{section}
\allowdisplaybreaks[2]
\setlist[enumerate]{label=\textup{(\roman*)},leftmargin=*}
\title[Rank Additivity and Comon's Conjecture for Real Hankel Tensors]
{Rank Additivity and Comon's Conjecture for Real\\  Hankel Tensors}
\author{Jianze Li}
\address{School of Science, Sun Yat-sen University, Shenzhen, Guangdong, China}
\email{lijianze@mail.sysu.edu.cn}
\date{9 November 2026}

\begin{document}
\begin{abstract}
Let $\set{H}_m(n)$ denote the set of real Hankel tensors of order $m$ and dimension $n$, and let $\tens{H}\in\set{H}_m(n)$, where $m\ge3$ and $n\ge1$.
Set $M=(m-1)(n-1)+1$.
First, we prove that the CP rank $\rk(\tens{H})\le M$ and that CP rank is additive under direct sums of $\tens{H}$ with any real tensor $\tens{T}$ of the same order.
Second, we prove that partially symmetric rank $p_j(\tens{H})=\rk(\tens{H})$ for $1\le j\le m-1$ and that, if $\rk(\tens{H})\le M-1$, these ranks also equal its symmetric rank $s(\tens{H})$.
Third, we construct a fourth-order Hankel tensor $\tens{H}_*\in\set{H}_4(4)$ that attains this CP rank bound and satisfies
$\rk(\tens{H}_*)=10$ and $s(\tens{H}_*)=11$.
Moreover, for the direct sum $\tens{H}_*^{(k)}=\bigoplus_{a=1}^k\tens{H}_*$, we further prove that
$\rk(\tens{H}_*^{(k)})=10k$ and $s(\tens{H}_*^{(k)})=11k$.
This provides a low-dimensional counterexample to Comon’s conjecture over the real numbers and further shows that the Comon rank gap is unbounded among real symmetric tensors of order four. 
\end{abstract}
\maketitle

\section{Introduction}\label{sec:introduc}
\subsection{Tensor ranks and direct sums} 
The CP rank \cite{hitchcock1927expression,kruskal1977three,kolda2009tensor,landsberg2012tensors},
partially symmetric rank \cite{gesmundo2019partially,bernardi2013general,wang2023lower}, and symmetric rank \cite{comon2008symmetric,landsberg2012tensors}, which is
equivalently the \emph{Waring rank} of the associated homogeneous polynomial, are fundamental measures
of the complexity of tensor decompositions under different symmetry constraints. 
In this paper, unless otherwise stated, all vector spaces are finite-dimensional over $\RR$, and all tensor ranks are computed over $\RR$. 
In all rank definitions below, the minimum is taken over
$N\in\mathbb Z_{\ge0}$, and an empty sum represents the zero tensor. 
Let $m\ge3$ and let $\set{A}_1,\ldots,\set{A}_m$ be vector spaces.
For a tensor $\tens{T}\in\bigotimes_{i=1}^m\set{A}_i$, its \emph{CP rank} is
\[
\rk(\tens{T})\eqdef
\min\left\{N:\tens{T}=\sum_{\nu=1}^{N}
\vect{a}_{\nu,1}\otimes\cdots\otimes\vect{a}_{\nu,m},\
\vect{a}_{\nu,i}\in\set{A}_i\right\}.
\]
For a vector space $\set{A}$, let $\Sym^m(\set{A})$ denote the subspace
of $\set{A}^{\otimes m}$ invariant under all permutations of the factors.
The homogeneous polynomial associated with $\tens{S}\in\Sym^m(\set{A})$
is $\alpha\mapsto\langle\tens{S},\alpha^{\otimes m}\rangle$ on $\set{A}^*$.
For a symmetric tensor $\tens{S}\in\Sym^m(\set{A})$ and
$1\le k\le m$, define its \emph{partially symmetric ranks} and \emph{symmetric rank} by
\begin{align*}
p_k(\tens{S})
\eqdef
\min\left\{N:\tens{S}=\sum_{\nu=1}^{N}
\gamma_{\nu}\vect{a}_{\nu,1}\otimes\cdots\otimes
\vect{a}_{\nu,m-k}\otimes\vect{v}_{\nu}^{\otimes k}
\right\},\ 
s(\tens{S})
\eqdef
\min\left\{N:
\tens{S}=\sum_{\nu=1}^{N}
\gamma_{\nu}\vect{v}_{\nu}^{\otimes m}
\right\}.
\end{align*} 
Here $\gamma_{\nu}\in\RR$, $\vect{v}_{\nu}\in\set{A}$, and
$\vect{a}_{\nu,j}\in\set{A}$ for $1\le j\le m-k$.
When $k=m$, the factors preceding
$\vect{v}_{\nu}^{\otimes k}$ are omitted.
Since $\tens{S}$ is symmetric, the $k$ equal factors in the definition
of $p_k(\tens{S})$ may be placed in any fixed set of $k$ positions,
not necessarily the last $k$ positions, without changing its value.
It follows directly from the definitions that
\begin{equation}\label{eq:rank-chain}
\rk(\tens{S})=p_1(\tens{S})\le p_2(\tens{S})
\le\cdots\le p_m(\tens{S})=s(\tens{S}).
\end{equation} 

For tensors $\tens{T}_1\in\bigotimes_{i=1}^m\set{A}_i$ and
$\tens{T}_2\in\bigotimes_{i=1}^m\set{B}_i$,
the \emph{direct sum tensor} $\tens{T}_1\oplus\tens{T}_2$ belongs to
$\bigotimes_{i=1}^m(\set{A}_i\oplus\set{B}_i)$.
Its two diagonal blocks are $\tens{T}_1$ and $\tens{T}_2$, and all
mixed blocks vanish.
The above ranks are all subadditive under direct sums whenever
the corresponding definitions apply.
We define the \emph{Strassen CP rank gap}, \emph{Strassen symmetric rank gap} \cite{li2026hankeladditivity,li2026directsumgap},
and \emph{Comon rank gap} \cite{li2026comonrankgap} by
\begin{align*}
\Delta_r(\tens{T}_1,\tens{T}_2)
&\eqdef \rk(\tens{T}_1)+\rk(\tens{T}_2)
-\rk(\tens{T}_1\oplus\tens{T}_2),\\
\Delta_s(\tens{S}_1,\tens{S}_2)
&\eqdef s(\tens{S}_1)+s(\tens{S}_2)
-s(\tens{S}_1\oplus\tens{S}_2),\ \
\Delta_c(\tens{S})\eqdef s(\tens{S})-\rk(\tens{S}).
\end{align*}
The latter two definitions concern symmetric tensors.
A Strassen rank gap vanishes exactly when the corresponding rank
is additive for that pair of tensors, 
and the Comon rank gap vanishes exactly when the CP rank and symmetric rank coincide.

\emph{Strassen's direct sum conjecture} \cite{strassen1973division},
originally formulated for third-order tensors, asks whether
$\Delta_r(\tens{T}_1,\tens{T}_2)=0$ for every pair of such tensors.
Various sufficient conditions for CP rank additivity have been established
\cite{jaja1986validity,buczynski2020additivity,rupniewski2024additivity}.
However, Shitov disproved the general conjecture over every infinite
field in 2019 \cite{shitov2019directsum}. 
The \emph{Strassen's symmetric direct sum conjecture} asks whether
$\Delta_s(\tens{S}_1,\tens{S}_2)=0$ for every pair of symmetric tensors
of the same order \cite{carlini2015symmetricstrassen}.
Sufficient conditions have been obtained using rank lower bounds
and Waring loci \cite{teitler2015sufficient,carlini2017waringloci}.
Shitov's preprint \cite[Theorem~3.9]{shitov2023higher} reports
counterexamples to symmetric rank additivity over every infinite field
of characteristic different from two and three. 

\emph{Comon's conjecture}
\cite{comon2004canonical,comon2008symmetric}
asks whether $\Delta_c(\tens{S})=0$ for every symmetric tensor $\tens{S}$.
The equality has been established in several special cases over
$\RR$ or $\CC$, under additional rank or decomposition assumptions
\cite{comon2008symmetric,zhang2016comon,friedland2016remarks,seigal2020ranks}.
Related equalities among partially symmetric ranks are studied in
\cite{gesmundo2019partially}. 
A proposed complex
counterexample \cite{shitov2018counterexample} was subsequently shown to have an
error \cite{draisma2024erratum}; new counterexamples over the complex numbers
were recently reported in \cite{shitovhigher,Lovitz26:arxiv}. 
Over $\RR$, Shitov's preprint \cite{shitov2020comon} reports a
fourth-order counterexample, and Wang and Seigal \cite{wang2023lower}
establish a sixth-order counterexample.
Later preprints report third-order counterexamples over both
$\RR$ and $\CC$ in sufficiently large dimension \cite{shitovhigher}
and, subsequently, in dimension $27$ \cite{Lovitz26:arxiv}.
For a survey of earlier developments on Comon's and Strassen's conjectures,
see \cite{casarotti2020comon}. 

\subsection{Hankel tensors and related work}

A tensor $\tens{H}\in(\RR^n)^{\otimes m}$ is \emph{Hankel}
\cite{luque2003hankel,qi2015hankel,nie2019hankel} if there exist
real numbers $h_0,\ldots,h_{m(n-1)}$ such that
$H_{i_1\ldots i_m}=h_{i_1+\cdots+i_m-m}$ for
$1\le i_1,\ldots,i_m\le n$.
Every Hankel tensor is symmetric, and the sequence $(h_0,\ldots,h_{m(n-1)})$ is called the
\emph{generating sequence} of $\tens{H}$.
Hankel tensors are natural higher-order generalizations of the Hankel matrices
\cite{iohvidov1982hankel,usevich2014variable},
and their Vandermonde decompositions and positivity properties
over $\RR$ have been studied in
\cite{qi2015hankel,xu2016hankel,ding2017inheritance}.

The CP rank and symmetric rank of Hankel tensors have been studied
over $\CC$.
Nie and Ye \cite{nie2019hankel} proved that these ranks coincide
for generic complex Hankel tensors of even order or order three,
and formulated Comon's conjecture for all complex Hankel tensors
\cite[Conjecture~6.3]{nie2019hankel}.
Colbrook's repository manuscript \cite{colbrook2026generic}
presents a proof extending the generic equality to the remaining odd orders.
Holden's repository manuscript \cite{holden2026exact}
presents an exact rank formula and a proof of equality of
CP rank and symmetric rank for every complex Hankel tensor.
The CP rank of any complex Hankel tensor of order $m\ge3$
is proved to be additive under direct sum with any complex tensor of the
same order \cite{li2026hankeladditivity}. 

To our knowledge, there are relatively few studies on rank additivity and Comon's conjecture 
for real Hankel tensors.
Although related results have been obtained over $\CC$, they do not
directly imply the corresponding statements over $\RR$.
Indeed, a real tensor may have a shorter decomposition over $\CC$,
and apolarity arguments over $\RR$ must account for real roots and
irreducible quadratic factors.
Therefore, further research is needed to better understand
the ranks of Hankel tensors over $\RR$.

\subsection{Contributions}

In this paper, we study rank additivity and Comon's conjecture for
Hankel tensors over $\RR$.
We denote by $\set{H}_m(n)\subseteq (\RR^n)^{\otimes m}$ the space of real Hankel tensors
of order $m\ge3$ and dimension $n\ge1$, and set
\begin{equation}\label{eq:M_value}
M\eqdef(m-1)(n-1)+1.
\end{equation} 
Our main results consist of the following three parts. 
First, we prove that every $\tens{H}\in\set{H}_m(n)$ satisfies
$\rk(\tens{H})\le M$, and its CP rank is additive under direct sums
with any real tensor of the same order
\textup{(\cref{thm:main-additivity})}. 
Second, we prove that $p_j(\tens{H})=\rk(\tens{H})$ for every
$1\le j\le m-1$.
Moreover, when $\rk(\tens{H})\le M-1$, the symmetric rank also
equals the CP rank
\textup{(\cref{thm:main-comon})}. 
These results also imply symmetric rank additivity for the
direct sum of a real Hankel tensor of CP rank at most $M-1$
and any real symmetric tensor satisfying Comon rank equality 
\textup{(\cref{cor:symmetric-partner})}. 
Third, we construct a fourth-order Hankel tensor
$\tens{H}_*\in\set{H}_4(4)$ with entries
$(\tens{H}_*)_{ijkl}=\delta_{i+j+k+l,10}$ and prove that
$\rk(\tens{H}_*)=10$ and $s(\tens{H}_*)=11$. 
Furthermore, for the repeated direct sum
$\tens{H}_*^{(k)}=\bigoplus_{a=1}^k\tens{H}_*$, we prove the exact
formulas $\rk(\tens{H}_*^{(k)})=10k$ and $s(\tens{H}_*^{(k)})=11k$
for every $k\ge1$
\textup{(\cref{thm:main-counterexample})}. 
It follows that $\Delta_c(\tens{H}_*^{(k)})=k$.
Thus, the Comon rank gap is unbounded among real symmetric
tensors of order four. 
In particular, taking $k=2$ gives a positive answer to the real fourth-order analogue of \cite[Problem 23]{shitov2018counterexample}:
there exists a real symmetric tensor with Comon rank gap equal to two. 

The rest of the paper is organized as follows.
\Cref{sec:main-results} states the main results.
\Cref{sec:complete-moment} develops moment representations and
computes the ranks of complete moment tensors.
\Cref{sec:apolar-tools} establishes apolar lifting and proper-support
reductions, and \cref{sec:decomposition-methods} develops decomposition
blocks and support bounds.
\Cref{sec:restriction} proves the symmetric rank restriction estimate
and CP rank additivity, while \cref{sec:critical} establishes the
rank equalities.
\Cref{sec:quartic-counterexample} determines the exact ranks
of $\tens{H}_*$, and
\cref{sec:quartic-direct-sum-amplification} treats its repeated
direct sums and proves the unbounded-gap result.
Finally, \cref{sec:outlook} summarizes the results and presents open problems.

\section{Main results}\label{sec:main-results}

This section presents our main results in three parts.
We denote by $\set{H}_m(n)\subseteq (\RR^n)^{\otimes m}$
the space of real Hankel tensors of order $m$ and dimension $n$.

\subsection{CP rank additivity}

\begin{theorem}[Rank bound and direct-sum additivity]
\label{thm:main-additivity}
Let $\tens{H}\in\set{H}_m(n)$ with $m\ge3$ and $n\ge1$, and let $M$ be as in
\eqref{eq:M_value}. Then
$\rk(\tens{H})\le M.$ 
Moreover, for any real tensor
$\tens{T}\in\set{B}_1\otimes\cdots\otimes\set{B}_m$, we have $\Delta_r(\tens{H},\tens{T})=0$, that is, 
\begin{equation}\label{eq:main-additivity}
\rk(\tens{H}\oplus\tens{T})
=\rk(\tens{H})+\rk(\tens{T}).
\end{equation}
\end{theorem}

\subsection{Rank equalities}

\begin{theorem}[Rank equalities for real Hankel tensors]
\label{thm:main-comon}
Let $\tens{H}\in\set{H}_m(n)$ with $m\ge3$ and $n\ge1$, and let $M$ be as in
\eqref{eq:M_value}. Then
\[
\rk(\tens{H})=p_1(\tens{H})
=\cdots=p_{m-1}(\tens{H}).
\]
If $\rk(\tens{H})\le M-1$, then we further have that 
\[
\rk(\tens{H})=p_1(\tens{H})
=\cdots=p_{m-1}(\tens{H})=p_m(\tens{H})=s(\tens{H}).
\]
\end{theorem}

\begin{corollary}[A symmetric partner satisfying Comon rank equality]
\label{cor:symmetric-partner}
Let $\tens{H}\in\set{H}_m(n)$ with $m\ge3$ and $n\ge1$, and let $M$ be as in
\eqref{eq:M_value}.
Let $\tens{S}\in\Sym^m(\set{A})$ satisfy
$s(\tens{S})=\rk(\tens{S})$. If
$\rk(\tens{H})\le M-1$, then
\begin{align*}
\rk(\tens{H}\oplus\tens{S})
=s(\tens{H}\oplus\tens{S}) 
=\rk(\tens{H})+\rk(\tens{S})
=s(\tens{H})+s(\tens{S}).
\end{align*}
\end{corollary}

\subsection{The maximal rank obstruction and its amplification}

Let $\tens{H}_*\in\set{H}_4(4)$ be the Hankel tensor defined, for
$1\le i,j,k,l\le4$, by
\begin{equation}\label{eq:main-counterexample-tensor}
(\tens{H}_*)_{ijkl}\eqdef\delta_{i+j+k+l,10}.
\end{equation}
For each integer $k\ge1$, set
\begin{equation}\label{eq:main-counterexample-tensor-k}
\tens{H}_*^{(k)}\eqdef\bigoplus_{a=1}^k\tens{H}_*
\in\Sym^4\left(\bigoplus_{a=1}^k\RR^4\right).
\end{equation}
Thus $\tens{H}_*^{(k)}$ is a fourth-order real symmetric tensor
of dimension $4k$.

\begin{theorem}[Exact maximal obstruction and unbounded gaps]
\label{thm:main-counterexample}
Let $\tens{H}_*\in\set{H}_4(4)$ be as in
\eqref{eq:main-counterexample-tensor}. Then
\[
\rk(\tens{H}_*)=p_2(\tens{H}_*)=p_3(\tens{H}_*)=10,
\ \  s(\tens{H}_*)=11.
\]
For every $k\ge1$, its repeated direct sum in
\eqref{eq:main-counterexample-tensor-k} satisfies
\begin{equation}\label{eq:main-direct-sum-gap}
\begin{aligned}
\rk(\tens{H}_*^{(k)})=p_2(\tens{H}_*^{(k)})=p_3(\tens{H}_*^{(k)})=10k,\ \
s(\tens{H}_*^{(k)})=11k.
\end{aligned}
\end{equation}
Consequently $s(\tens{H}_*^{(k)})-\rk(\tens{H}_*^{(k)})=k$.
Thus, the Comon rank gap is unbounded among fourth-order
real symmetric tensors, and every positive integer occurs
as such a gap.
\end{theorem}

\section{Moment representations and complete moment tensors}
\label{sec:complete-moment}

\subsection{Shortest decompositions and supports}

The \emph{length} of a decomposition is its number of nonzero terms.
For a fixed decomposition, a \emph{subdecomposition} means a sublist of its
terms, and a \emph{subsum} is the tensor represented by that sublist. 
For a tensor $\tens{T}\in\bigotimes_{j=1}^m\set{A}_j$, its
\emph{mode-$j$ support} $\supp_j(\tens{T})\subseteq\set{A}_j$ is the image of the
mode-$j$ contraction map
$\bigotimes_{i\ne j}\set{A}_i^*\to\set{A}_j$.
The mode-$j$ \emph{flattening} is the matrix of this contraction map
in chosen bases; the \emph{first flattening} corresponds to $j=1$.
The tensor is called \emph{concise} if
$\supp_j(\tens{T})=\set{A}_j$ for every $j$.
We use the following standard facts about shortest decompositions,
stated here for tensors of arbitrary order $m\ge3$; see
\cite[Lemma 2.8]{buczynski2020additivity} and \cite[Lemma~2.1]{li2026directsumgap} for the third-order case.

\begin{lemma}[Minimum decompositions]
\label{lem:shortest-subsum-support}
Let $m\ge3$ and let $\tens{T}\in\bigotimes_{j=1}^m\set{A}_j$ have a decomposition
$\tens{T}=\sum_{\nu=1}^{N}
\vect{a}_{\nu,1}\otimes\cdots\otimes\vect{a}_{\nu,m},$ 
where $N=\rk(\tens{T})$ and $\vect{a}_{\nu,j}\in\set{A}_j$.
For every nonempty $\set{I}\subseteq\{1,\ldots,N\}$, set
$\tens{T}_{\set{I}}=\sum_{\nu\in \set{I}}
\vect{a}_{\nu,1}\otimes\cdots\otimes\vect{a}_{\nu,m}$. Then $\rk(\tens{T}_{\set{I}})=|\set{I}|$ and 
$$\supp_j(\tens{T}_{\set{I}})
=\Span\{\vect{a}_{\nu,j}:\nu\in \set{I}\},\ \ 1\le j\le m.$$ 
Similarly, let $\tens{S}\in\Sym^m(\set{A})$ have a decomposition
$\tens{S}=\sum_{\nu=1}^{N}\gamma_{\nu}\vect{v}_{\nu}^{\otimes m},$ 
where $N=s(\tens{S})$, $\gamma_{\nu}\in\RR\setminus\{0\}$, and $\vect{v}_{\nu}\in\set{A}$.
For every nonempty $\set{I}\subseteq\{1,\ldots,N\}$, the subsum
$\tens{S}_{\set{I}}=\sum_{\nu\in \set{I}}\gamma_{\nu}\vect{v}_{\nu}^{\otimes m}$ satisfies $s(\tens{S}_{\set{I}})=|\set{I}|$ and 
$$\supp_j(\tens{S}_{\set{I}})=\Span\{\vect{v}_{\nu}:\nu\in \set{I}\},\ \ 1\le j\le m.$$ 
\end{lemma}

\subsection{Moment representation and apolar notation}

Let $\set{P}_j\eqdef\RR[x,y]_j$ be the space of homogeneous binary forms
of degree $j$ for $j\ge0$, and set $\set{P}_j\eqdef0$ for $j<0$.
For $j\ge0$, we identify $\set{P}_j$ with $\RR[t]_{\le j}$, the space of real univariate polynomials of degree at most $j$, via $f(x,y)\mapsto f(1,t)$. 
Products and functionals are interpreted under this identification.
For polynomial subspaces $\set{U}$ and $\set{V}$, write
$\set{U}\set{V}\eqdef
\Span\{uv:u\in\set{U},\,v\in\set{V}\}$.
For a fixed polynomial $g$, write
$g\set{U}\eqdef\{gu:u\in\set{U}\}$.
Products of more than two polynomial subspaces are defined similarly.
For a nonzero polynomial subspace $\set{K}$, its
\emph{greatest common divisor} (gcd) is the gcd of the polynomials
in any basis of $\set{K}$, defined up to a nonzero scalar.
For a subspace of homogeneous binary forms, the gcd is computed in
$\RR[x,y]$ before dehomogenization, so factors at infinity are retained.
If a nonzero polynomial $h$ divides every element of $\set{K}$, write
$\set{K}/h\eqdef\{f/h:f\in\set{K}\}$.
The notation $[t^j]f$ denotes the coefficient of $t^j$ in $f$.

Let $q\in\RR[t]$ be monic and nonconstant, where $\RR[t]$ is
the ring of real univariate polynomials.
Define the \emph{quotient algebra} $\set{A}_q\eqdef\RR[t]/(q)$,
where $(q)\eqdef\{qf:f\in\RR[t]\}$ is the ideal generated by $q$.
Two polynomials represent the same element of $\set{A}_q$
if and only if their difference is divisible by $q$.
The notation $[f]_q$ denotes the unique remainder of $f$ modulo $q$,
which has degree less than $\deg q$.
We identify this representative with its residue class in $\set{A}_q$.
Addition and multiplication in $\set{A}_q$ are computed modulo $q$,
and $\dim\set{A}_q=\deg q$.

For a vector space $\set{V}$ and $\set{S}\subseteq\set{V}^*$, write
$\set{S}^\perp\eqdef\{v\in\set{V}:\alpha(v)=0\text{ for every }\alpha\in\set{S}\}$.
For a subspace $\set{U}\subseteq\set{V}$, its \emph{annihilator} in
$\set{V}^*$ is defined by
$\set{U}^\perp\eqdef
\{\ell\in\set{V}^*:\ell(u)=0\text{ for every }u\in\set{U}\}$.
Let $n\ge1$ and $D=m(n-1)$.
Use the ordered basis
$(x^{n-1},x^{n-2}y,\ldots,y^{n-1})$ of $\set{P}_{n-1}$ and its dual basis
to identify $(\set{P}_{n-1}^*)^{\otimes m}$ with $(\RR^n)^{\otimes m}$.
If $(h_0,\ldots,h_D)$ is the generating sequence of $\tens{H}$,
define the \emph{moment functional} $L\in \set{P}_D^*$ by $L(x^{D-j}y^j)\eqdef h_j$.
Then the indexing convention in \cref{sec:introduc} gives the
moment representation
\begin{equation*}
\tens{H}(f_1,\ldots,f_m)=L(f_1\cdots f_m),
\end{equation*}
where $f_i\in \set{P}_{n-1}$ and $L\in \set{P}_D^*$.
For a linear subspace $\set{K}\subseteq \set{P}_{n-1}$, write
$\tens{H}_{\set{K}}\eqdef\tens{H}|_{\set{K}\times \set{P}_{n-1}^{m-1}}$.
The tensor is concise if its support in every mode equals
$\set{P}_{n-1}^*$.

For a nonzero real univariate polynomial or homogeneous binary form $g$, let $\tau(g)$ be the number of distinct irreducible factors over $\RR$, without counting multiplicities. For a homogeneous binary form, $\deg g$ and $\tau(g)$ are computed from the homogeneous form, so any factor $x$ is counted even though it disappears after setting $x=1$. 
We define the \emph{apolar cost} of $g$ as 
\begin{equation}\label{eq:gamma}
\Gamma_m(g)\eqdef(m-1)\deg g-(m-2)\tau(g)
=\deg g+(m-2)(\deg g-\tau(g)).
\end{equation}
Set $\tau(1)=\Gamma_m(1)=0$. 
For a nonconstant univariate polynomial $g$, the expression
$\Gamma_m(g)$ is the specialization of the multilinear
Alder--Strassen lower bound \cite{alder1981algorithmic} to $(m-1)$-fold multiplication
in $\RR[t]/(g)$ \cite[Theorem~19]{blaser2023multilinear}. 
A nonzero homogeneous binary form $g$ is \emph{apolar} to $\tens{H}$ if $L(gp)=0$ for every $p\in \set{P}_{D-\deg g}$. 
The \emph{apolar ideal} of $L$ is the homogeneous ideal whose degree-$r$
part is $\{g\in\set{P}_r:L(gp)=0\text{ for every }p\in\set{P}_{D-r}\}$
for $r\ge0$.
The following properties follow directly from factorization over $\RR$. 

\begin{lemma}[Elementary properties of the apolar cost]
\label{lem:apolar-cost-arithmetic}
Let $m\ge3$ and let $g,h$ be nonzero real binary forms. Then $\deg g\le\Gamma_m(g)$ and 
\begin{align*}
g\mid h\ \Longrightarrow\ \Gamma_m(g)\le\Gamma_m(h);\ 
\Gamma_m(\lcm(g,h))\le\Gamma_m(g)+\Gamma_m(h);\ 
\Gamma_m(gh)\le\Gamma_m(g)+(m-1)\deg h.
\end{align*}
\end{lemma}

\subsection{Universal rank bounds}
\label{subsec:universal-bounds}


\begin{proposition}[Universal CP and partially symmetric rank bounds]
\label{prop:partial-upper}
\label{cor:ordinary-upper}
Let $\tens{H}\in\set{H}_m(n)$ with $m\ge3$ and $n\ge1$,
let $L\in\set{P}_{m(n-1)}^*$ be its moment functional,
and let $M$ be as in \eqref{eq:M_value}. Then
\begin{equation*}
\rk(\tens{H})\le p_{m-1}(\tens{H})\le M.
\end{equation*}
\end{proposition}

\begin{proof}
Choose $M$ distinct real nodes $\theta_1,\ldots,\theta_M$, and let $\ell_i\in\RR[t]_{\le M-1}$ be the corresponding Lagrange polynomials. 
Under the identification $\set{P}_{n-1}\simeq\RR[t]_{\le n-1}$,
the product $f_2\cdots f_m$ has degree at most $M-1$.
Lagrange interpolation therefore gives
$f_2\cdots f_m
=\sum_{i=1}^{M}\left(\prod_{j=2}^{m}f_j(\theta_i)\right)\ell_i.$ 
Consequently
\[
\tens{H}(f_1,\ldots,f_m)
=\sum_{i=1}^{M}L(f_1\ell_i)\prod_{j=2}^{m}f_j(\theta_i).
\]
The same evaluation functional occurs in the last $m-1$ modes
of each term.
Thus this is a decomposition of the type defining
$p_{m-1}(\tens{H})$, with at most $M$ terms.
\end{proof}

\subsection{Exact ranks of complete moment tensors}
\label{subsec:complete-moment-ranks}

Let $m\ge3$, put $\set{A}_q\eqdef\RR[t]/(q)$, and let $\lambda:\set{A}_q\to\RR$
be a \emph{Frobenius functional}, meaning that the bilinear form
$(u,v)\mapsto\lambda(uv)$ is nondegenerate.  
The associated \emph{complete moment tensor} is defined as 
\begin{equation}\label{eq:def_Tm}
\tens{T}_{m,q,\lambda}(u_1,\ldots,u_m)
\eqdef \lambda(u_1\cdots u_m)
\ \ (u_i\in \set{A}_q). 
\end{equation}
Here completeness means that every input ranges over
the whole algebra.
The $(m-1)$-fold \emph{multiplication tensor} of $\set{A}_q$ is the
$m$th-order tensor representing the multilinear map
$(u_1,\ldots,u_{m-1})\mapsto u_1\cdots u_{m-1}$
from $\set{A}_q^{m-1}$ to $\set{A}_q$.
Under the Frobenius identification
$\set{A}_q\simeq \set{A}_q^*$, given by
$z\mapsto\lambda(z\,\cdot)$,
the $(m-1)$-fold multiplication tensor of $\set{A}_q$
becomes $\tens{T}_{m,q,\lambda}$.
In particular, their CP ranks coincide;
see \cite{blaser2023multilinear}. 
For $q=t^e$ and $\lambda=[t^{e-1}]$ on remainder
representatives, its complexification is precisely
the tensor $\tens{L}(e,m)$ studied in
\cite[Section~4]{gharahi2024persistent}.

\begin{lemma}[Primary real power cubature]
\label{prop:primary-cubature}
Let $m\ge3$, let $q=f^a$ be monic with $f\in\RR[t]$ irreducible and
$a\in\mathbb Z_{\ge1}$, put $\set{A}_q=\RR[t]/(q)$,
and let $\lambda:\set{A}_q\to\RR$ be Frobenius.  Put
$N=(m-1)(\deg q-1)+1$.  Then the homogeneous polynomial
$u\mapsto\lambda(u^m)$ is a signed sum of at most $N$ real $m$th
powers of linear forms on $\set{A}_q$.
\end{lemma}

\begin{proof}
Put $e=\deg q$. 
If $e=1$, then $\set{A}_q\simeq\RR$ and the assertion is immediate.
Assume henceforth that $e\ge2$. 
Let $\lambda_0$ extract the coefficient of $t^{e-1}$ from the
canonical representative of degree less than $e$.  Its multiplication
pairing is nondegenerate: if $0\ne u\in \set{A}_q$ has leading degree $j$,
then $\lambda_0(ut^{e-1-j})\ne0$.  
Consequently, there is a unique element $w\in \set{A}_q$ such that 
\begin{equation}\label{eq:frob-unit}
\lambda(u)=\lambda_0(wu)\ \ (u\in \set{A}_q).
\end{equation}
Since $\lambda$ is Frobenius, multiplication by $w$ is injective.
As $\set{A}_q$ is finite-dimensional, this map is invertible, so $w$
is a unit. 

Choose distinct real numbers $\theta_1,\ldots,\theta_N$ away from the
real roots of $q$ and let $h(t)=\prod_{i=1}^N(t-\theta_i)$.  Its class
in $\set{A}_q$ is a unit. The residue field $\set{A}_q/(f)$ is $\RR$ or $\CC$.
Choose $\varepsilon\in\{1,-1\}$ so that the residue of
$\varepsilon(wh)^{-1}$ has an $m$th root, and lift that root to a unit $z_0$.
Only a negative real residue with even $m$ requires $\varepsilon=-1$.
Put $\eta\eqdef\varepsilon(wh)^{-1}z_0^{-m}-1$.
Then $\eta^a=0$, so the finite binomial expression
$z\eqdef z_0\sum_{j=0}^{a-1}\binom{1/m}{j}\eta^j$
satisfies $z^m=\varepsilon(wh)^{-1}$.
Equation \eqref{eq:frob-unit} now gives
$\lambda((zu)^m)=\varepsilon\lambda_0(h^{-1}u^m)$. 
For a polynomial $P$ with $\deg P\le e+N-2$, put
$Q_P\eqdef(P-h[h^{-1}P]_q)/q$.
Then $\deg[h^{-1}P]_q<e$ and $\deg Q_P<N$.
Comparing the coefficient of $t^{e+N-1}$ in $P=qQ_P+h[h^{-1}P]_q$ gives
$[t^{e-1}][h^{-1}P]_q=-[t^{N-1}]Q_P$.
At each node we have $Q_P(\theta_i)=P(\theta_i)/q(\theta_i)$.
Taking the leading coefficient in the Lagrange formula for $Q_P$ yields
\[
\lambda_0(h^{-1}P)
=-\sum_{i=1}^N\frac{P(\theta_i)}{q(\theta_i)h'(\theta_i)}.
\]
For $u\in\set{A}_q$, apply this identity to $[z^{-1}u]_q^m$.
Since $\deg([z^{-1}u]_q^m)\le m(e-1)=e+N-2$, it follows that
\[
\lambda(u^m)
=\varepsilon\lambda_0(h^{-1}[z^{-1}u]_q^m)
=-\varepsilon\sum_{i=1}^N
\frac{\bigl([z^{-1}u]_q(\theta_i)\bigr)^m}{q(\theta_i)h'(\theta_i)}.
\]
Each map $u\mapsto([z^{-1}u]_q)(\theta_i)$ is a real linear form on $\set{A}_q$.
Absorbing the $m$th root of the absolute value of each coefficient
into the corresponding linear form gives the required signed power decomposition.
\end{proof}
For $m=3$, the Frobenius pairing identifies the moment tensor with
the multiplication tensor of $\set{A}_q$, and the lower bound below is defined in \eqref{eq:gamma}.
The following proof also covers higher orders and supplies signed
real symmetric decompositions attaining the bound.

\begin{proposition}[Exact ranks of a complete moment tensor]
\label{thm:complete-rank}
Let $q\in\RR[t]$ be monic of degree $e\ge1$, and let
$\lambda$ be a Frobenius functional on $\set{A}_q=\RR[t]/(q)$.
For the tensor $\tens{T}_{m,q,\lambda}$ defined in
\eqref{eq:def_Tm}, we have
\[
\rk(\tens{T}_{m,q,\lambda})
=(m-1)e-(m-2)\tau(q)
=\Gamma_m(q)
\ \  (m\ge3).
\]
The same number equals every
$p_j(\tens{T}_{m,q,\lambda})$, $1\le j\le m$,
and hence $s(\tens{T}_{m,q,\lambda})$.
\end{proposition}

\begin{proof}
Write
$q=\prod_{j=1}^{\tau(q)}q_j$
as a product of pairwise coprime real \emph{primary factors},
each a positive integer power of a monic irreducible polynomial.
Put $\set{A}_{q_j}=\RR[t]/(q_j)$, and let
$\iota_j:\set{A}_{q_j}\to \set{A}_q$ be the canonical embedding of the
$j$th factor under the Chinese remainder isomorphism.
Define
$\lambda_j=\lambda\circ\iota_j.$ 
Each $\lambda_j$ is Frobenius, since the multiplication
pairing is block diagonal. Under the same Chinese
remainder isomorphism in every mode, we have
$\tens{T}_{m,q,\lambda}
\cong
\bigoplus_{j=1}^{\tau(q)}
\tens{T}_{m,q_j,\lambda_j}.$ 
Applying \cref{prop:primary-cubature} to every primary
summand and concatenating the decompositions gives
\begin{align*}
s(\tens{T}_{m,q,\lambda})
\le
\sum_{j=1}^{\tau(q)}
s(\tens{T}_{m,q_j,\lambda_j}) 
\le
\sum_{j=1}^{\tau(q)}
\bigl((m-1)(\deg q_j-1)+1\bigr)
=(m-1)e-(m-2)\tau(q)
=\Gamma_m(q).
\end{align*}

For the reverse inequality, the Frobenius
identification makes $\tens{T}_{m,q,\lambda}$ equivalent to the
$(m-1)$-fold multiplication tensor of $\set{A}_q$.
Since $\dim\set{A}_q=e$ and $\set{A}_q$ has $\tau(q)$ maximal ideals,
the multilinear Alder--Strassen bound
\cite[Theorem~19]{blaser2023multilinear} gives
$$\rk(\tens{T}_{m,q,\lambda})\ge(m-1)e-(m-2)\tau(q)=\Gamma_m(q).$$
Combining this with the signed symmetric upper bound
and the rank chain
\eqref{eq:rank-chain} proves all asserted equalities.
\end{proof}

\begin{remark}
For $q=t^e$ and $\lambda=[t^{e-1}]$ on remainder
representatives, the corresponding complex CP rank
was computed in \cite[Theorem~3]{gharahi2024persistent}.
The proof above additionally provides matching signed
real symmetric decompositions.
\end{remark}

\section{Apolar lifting and proper-support reductions}
\label{sec:apolar-tools}

\subsection{Frobenius reduction and apolar quotients}

For linear maps $\phi_i:\set{B}_i\to\set{A}_i$ and an
$m$-linear form $\tens{T}$ on
$\set{A}_1\times\cdots\times\set{A}_m$, its \emph{pullback} is the
$m$-linear form on $\set{B}_1\times\cdots\times\set{B}_m$
given by
$$(\vect{v}_1,\ldots,\vect{v}_m)\mapsto
\tens{T}(\phi_1(\vect{v}_1),\ldots,\phi_m(\vect{v}_m)).$$
This is the standard pullback construction for covariant tensors, applied separately in each input;
see \cite[Chapter~12]{lee2012smooth}. 
The Frobenius reduction below is the univariate form of the
standard quotient construction associated with a Hankel
functional; see \cite{mourrain2018moments}.
We include a short proof, together with the elementary
rank-preserving pullback argument needed here.

\begin{lemma}[Frobenius reduction and rank-preserving pullback]
\label{lem:frobenius-reduction}
Let $m\ge3$, let $q\in\RR[t]$ be nonconstant, let
$\set{A}_q=\RR[t]/(q)$, and let $0\ne\mu\in \set{A}_q^*$.
There is a monic nonconstant divisor $g\mid q$ such that $\mu$ factors through $\set{A}_g=\RR[t]/(g)$ and the induced functional $\mu_g$ is Frobenius.  
If a finite-dimensional real polynomial window $\set{V}\subseteq\RR[t]$ maps onto
$\set{A}_g$ under reduction modulo $g$, let
$\pi:\set{V}\to\set{A}_g$ denote this map.
Then the pullback
$$(f_1,\ldots,f_m)\mapsto
\mu_g(\pi(f_1)\cdots\pi(f_m)),$$
viewed as a tensor in $(\set{V}^*)^{\otimes m}$,
has the same CP rank, partially symmetric ranks,
and real symmetric rank as the complete moment tensor
\[
\tens{T}_{m,g,\mu_g}(u_1,\ldots,u_m)
=\mu_g(u_1\cdots u_m),
\ \  u_i\in \set{A}_g. 
\]
\end{lemma}

\begin{proof}
Consider the \emph{multiplication radical}
$\set{J}_\mu\eqdef\{u\in \set{A}_q:\mu(u\set{A}_q)=0\}$.
It is an ideal contained in $\ker\mu$, and it is proper since
$\mu\ne0$.  Hence $\set{J}_\mu=(g)/(q)$ for a monic nonconstant
divisor $g\mid q$.  The induced functional on
$\set{A}_q/\set{J}_\mu\simeq\RR[t]/(g)$ is Frobenius by the
definition of the radical. 
Choose a linear section
$\sigma:\set{A}_g\to\set{V}$ with
$\pi\sigma=\mathrm{id}_{\set{A}_g}$.
Pullback along $\pi$ gives the window tensor, and pullback along
$\sigma$ recovers the complete tensor.  Rank cannot increase
under these linear maps.  Since the same map is used in every
mode, powers remain powers and equal tail factors remain equal.
Thus the two pullbacks give opposite inequalities for all the
asserted ranks.
\end{proof}

\begin{lemma}[Apolar identities and quotient moment tensors]
\label{lem:apolar-quotient}
Let $\tens{H}\in\set{H}_m(n)$ be a real Hankel tensor of order
$m\ge3$ and dimension $n\ge1$, with associated moment functional
$L\in\set{P}_D^*$, where $D=m(n-1)$.
Let $g$ be a nonzero binary form with $1\le\deg g\le D$
that is apolar to $\tens{H}$.
After a real change of binary coordinates, $\tens{H}$ is the
pullback, along the reduction map $\pi_g:\set{P}_{n-1}\to\set{A}_g$,
of a complete moment tensor with functional $\mu\in\set{A}_g^*$ on
$\set{A}_g=\RR[t]/(g)$, where $g$ now denotes its monic dehomogenization.
If $\mu\ne0$, its Frobenius reduction
is supported on a divisor $g_0\mid g$, and
\[
s(\tens{H})\le\Gamma_m(g_0)\le\Gamma_m(g).
\]
If $\mu\ne0$ and the reduction map $\set{P}_{n-1}\to\set{A}_{g_0}$
is surjective, the CP, partially symmetric, and symmetric ranks equal
the corresponding complete-tensor ranks. If $\mu=0$, then $\tens{H}=0$.
\end{lemma}

\begin{proof}
Choose real binary coordinates with $g(0,1)\ne0$ and normalize $g(1,t)$
to be monic. This common invertible change preserves all
ranks, degrees, and numbers of distinct irreducible factors.
Division by $g$ gives $L(f)=L([f]_g)$ for $f\in\RR[t]_{\le D}$
by apolarity. Thus $\mu([r]_g)\eqdef L(r)$ for $\deg r<\deg g$ defines
$\mu\in\set{A}_g^*$ with
\[
\tens{H}(v_1,\ldots,v_m)
=\mu\left(\prod_{j=1}^m[v_j]_g\right).
\]
The case $\mu=0$ is immediate. Otherwise \cref{lem:frobenius-reduction}
gives a divisor $g_0\mid g$ and a Frobenius functional $\mu_{g_0}$.
Pulling back a shortest signed power decomposition of
$\tens{T}_{m,g_0,\mu_{g_0}}$ gives
$$s(\tens{H})\le\Gamma_m(g_0)\le\Gamma_m(g)$$ by
\cref{thm:complete-rank,lem:apolar-cost-arithmetic}.
If the window map onto $\set{A}_{g_0}$ is surjective,
the rank-preserving part of \cref{lem:frobenius-reduction}
gives all the asserted equalities.
\end{proof}

\subsection{Binary saturation and apolar lifting}
\label{subsec:polynomial-tools}

The following saturation statement is a standard consequence
of the growth estimate for base-point-free binary linear systems
\cite[Lemma~9.8]{harris1992algebraic}.
We apply this estimate after complexification, which preserves dimensions.
The real gcd-one condition rules out a common complex projective zero,
since any nonreal common zero would also give its conjugate and hence
a common real quadratic factor.
Indeed, until $\set{U}\set{P}_d$ fills its ambient space,
this estimate increases its dimension by at least two at each
step, while the ambient dimension increases by one.
Its initial codimension $j+1-\dim\set{U}$ therefore gives the
stated threshold.

\begin{lemma}[Binary saturation]
\label{lem:binary-saturation}
Let $j,d\in\mathbb Z_{\ge0}$, let $0\ne \set{U}\subseteq \set{P}_j$, and assume that the
elements of $\set{U}$ have no nonconstant common divisor.  If $d\ge j+1-\dim\set{U}$, then
\begin{equation*}
\set{U}\set{P}_{d}=\set{P}_{j+d}.
\end{equation*}
\end{lemma}

\begin{proposition}[Apolar lifting with an explicit cost]
\label{lem:apolar-lifting}
Let $\tens{H}\in\set{H}_m(n)$ be a real Hankel tensor of order
$m\ge3$ and dimension $n\ge1$, with moment functional
$L\in\set{P}_{m(n-1)}^*$. 
Let $0\ne \set{K}\subseteq \set{P}_{n-1}$ have codimension $u$, let
$h=\gcd \set{K}$ have degree $a$, and let $0\ne f\in \set{P}_e$
with $e\in\mathbb Z_{\ge0}$.
Assume $$e+u-a\le(m-1)(n-1)\ \  \textrm{and}\ \  
L\bigl(\set{K}f \set{P}_{(m-1)(n-1)-e}\bigr)=0.$$   
Then $hf$ is apolar to $\tens{H}$ and
\begin{equation*}
s(\tens{H})\le\Gamma_m(f)+(m-1)a
\le\Gamma_m(f)+(m-1)u.
\end{equation*}
The same inequalities hold with $s$ replaced by $\rk$.
\end{proposition}

\begin{proof}
The inclusion $\set{K}\subseteq h\set{P}_{n-1-a}$ gives $a\le u$.  The space
$\set{K}/h\subseteq \set{P}_{n-1-a}$ has gcd one and codimension $u-a$.
The degree condition gives $(m-1)(n-1)-e\ge u-a$, so
\cref{lem:binary-saturation} gives
\[
\set{K}\set{P}_{(m-1)(n-1)-e}=h\set{P}_{m(n-1)-e-a}.
\]
Consequently the assumed annihilation is $L(hf\set{P}_{m(n-1)-e-a})=0$,
which is exactly the apolarity of $hf$.
If $hf$ is constant, this identity gives $L=0$ and the result is
immediate.  Otherwise \cref{lem:apolar-quotient} gives
$s(\tens{H})\le\Gamma_m(hf)$.  \cref{lem:apolar-cost-arithmetic} gives
\[
\Gamma_m(hf)\le\Gamma_m(f)+(m-1)a.
\]
The inequality $a\le u$ proves the remaining bound, and
$\rk(\tens{H})\le s(\tens{H})$ gives its CP rank version.
\end{proof}

\subsection{Common recurrences and nonconcise tensors}

For $q\ge2$, a tensor $\tens{Y}\in\set{B}\otimes(\RR^n)^{\otimes q}$
is \emph{Hankel in its tail modes} if contraction by every element
of $\set{B}^*$ gives a Hankel tensor in the remaining $q$ modes,
or a Hankel matrix when $q=2$.
A nonzero binary form $f$ is a \emph{common recurrence} of this family
if the moment functional of every such contraction annihilates
$f\set{P}_{q(n-1)-\deg f}$.

\begin{lemma}[Common recurrence of a proper Hankel family]
\label{lem:family-recurrence}
Let $m\ge3$, let $\set{B}$ be a real vector space,
and let $0\ne\tens{Y}\in \set{B}\otimes (\RR^n)^{\otimes(m-1)}$ be Hankel in its last
$m-1$ modes, and assume that their common support $\set{S}$ has dimension
$e<n$.  Then there is a real binary form $f$ of degree $e$ such that
the moment functional of every contraction of $\tens{Y}$ by an
element of $\set{B}^*$ annihilates $f\set{P}_{(m-1)(n-1)-e}$.
These are the coordinate moment functionals of the family.  Moreover
\begin{equation}\label{eq:Common recurrence}
\Gamma_m(f)\le (m-1)e-(m-2).
\end{equation}
\end{lemma}

\begin{proof}
Let $\set{R}\subseteq \set{P}_{n-1}$ be the annihilator of $\set{S}$, so that
$\dim \set{R}=n-e$, and put $f=\gcd \set{R}$, $e_0=\deg f$.
Every coordinate moment functional annihilates
$\set{R}\set{P}_{(m-2)(n-1)}$.  
Since $\set{R}\subseteq f\set{P}_{n-1-e_0}$ and $\dim \set{R}=n-e$,
we have $e_0\le e$.
The space 
$\set{R}/f=\{v/f:v\in \set{R}\}\subseteq \set{P}_{n-1-e_0}$ 
has gcd one and codimension $e-e_0$. 
Since $$e-e_0\le n-1\le(m-2)(n-1),$$ 
\cref{lem:binary-saturation} gives
$\set{R}\set{P}_{(m-2)(n-1)}=f\set{P}_{(m-1)(n-1)-e_0}$.
For every $v\in f\set{P}_{n-1-e_0}$, all coordinate moment functionals
vanish on $v\set{P}_{(m-2)(n-1)}$.  Thus contraction by $v$ in a tail
mode is zero, so $v\in \set{R}$.  
Hence $f\set{P}_{n-1-e_0}\subseteq \set{R}$.
Together with the inclusion $\set{R}\subseteq f\set{P}_{n-1-e_0}$,
this gives $\set{R}=f\set{P}_{n-1-e_0}$.
Comparing dimensions yields $n-e=n-e_0$, so $e_0=e$. 
Since a positive-degree real
form has at least one irreducible factor, formula \eqref{eq:gamma}
yields \eqref{eq:Common recurrence}.
\end{proof}

\begin{proposition}[Nonconcise Hankel tensors]
\label{prop:nonconcise-hankel}
Let $0\ne\tens{H}\in\set{H}_m(n)$ with $m\ge3$, and let
$L\in\set{P}_{m(n-1)}^*$ be its moment functional.
If its common mode support has dimension $1\le e<n$
(equivalently, its first flattening has rank $e$), then it is the pullback, in every mode, of a complete
Frobenius moment tensor on an algebra of dimension $e$.
Moreover,
\[
\rk(\tens{H})=p_{m-1}(\tens{H})=s(\tens{H}).
\]
\end{proposition}

\begin{proof}
Apply \cref{lem:family-recurrence}, with its order parameter taken to be
$m+1$, to $1\otimes\tens{H}$ with scalar coefficient space and $m$ tail modes.
It gives a real
binary form $h$ of degree $e$ that is apolar to $\tens{H}$.
After a common real change of binary coordinates,
\cref{lem:apolar-quotient} represents $\tens{H}$ as the pullback
of a complete Frobenius moment tensor
$\tens{T}_{m,g,\mu}$ on $\set{A}_g=\RR[t]/(g)$, where
$g\mid h$.  Since $\deg g\le e<n$, the reduction map
$\pi:\set{P}_{n-1}\to\set{A}_g$ is surjective.
The first flattening of this pullback has rank at most
$\dim\set{A}_g=\deg g$, so $e\le\deg g$.  Thus $\deg g=e$.
By \cref{lem:frobenius-reduction,thm:complete-rank},
$\tens{H}$ and $\tens{T}_{m,g,\mu}$ have the same CP, partially symmetric,
and symmetric ranks, and the asserted rank equalities follow.
\end{proof}

\section{Decomposition blocks and support bounds}
\label{sec:decomposition-methods}

\subsection{Detection components and exterior support}

The \emph{vector matroid} of a finite list of vectors has the linearly
independent sublists as its independent sets, and its rank is the
dimension of the span of the list.
For a nonempty list, this matroid is \emph{connected} if, across every
partition of the index set into two nonempty parts, the corresponding
spans have nonzero intersection.
Its \emph{connected components} are its maximal connected sublists;
zero columns form singleton components.
The following elementary lemma combines the connected-component
decomposition of a vector matroid
\cite[Corollary~4.2.9]{oxley2011matroid}
with rank-nullity.
We include a proof for completeness.

\begin{lemma}[Components of a linear detection map]
\label{lem:detector-components}
Let $\set{B},\set{V},\set{Z}$ be real vector spaces, and let
$\mathcal D:\set{V}\to \set{Z}$ be a real linear map.
Suppose $\tens{T}\in\set{B}\otimes\set{V}$ has a decomposition
$\tens{T}=\sum_{i=1}^N\vect{\alpha}_i\otimes\vect{t}_i$ for $\vect{\alpha}_i\in \set{B}$ and $\vect{t}_i\in \set{V}.$ 
Assume that
$(\operatorname{id}_{\set{B}}\otimes\mathcal D)\tens{T}=0$.
Partition the nonzero columns $\mathcal D(\vect{t}_i)$ into the connected components of their vector matroid over $\RR$, and treat each zero column as a singleton. 
For a component $\set{I}$, put
\[
\begin{aligned}
c_{\set{I}}\eqdef\dim\Span_{i\in \set{I}}\mathcal D(\vect{t}_i),\ \
k_{\set{I}}\eqdef\dim_{\RR}\Span_{i\in \set{I}}\vect{\alpha}_i.
\end{aligned}
\]
Then its subsum belongs to $\set{B}\otimes \set{V}_0$, where
$\set{V}_0\eqdef\ker\mathcal D$.
Moreover, with $k=\dim_{\RR}\Span_i\vect{\alpha}_i$, we have
\[
k_{\set{I}}+c_{\set{I}}\le |\set{I}|\ \ \text{and}\ \ \sum_{\set{I}}c_{\set{I}}\le N-k.
\]
If $\mathcal D$ instead varies polynomially with real parameters and
the assumed annihilation holds for every parameter choice, then
parameters can be chosen, after fixing the decomposition, so that
every component subsum belongs to $\set{B}\otimes\set{V}_0$, where
$\set{V}_0$ is now the common kernel, namely the intersection of the
kernels over all real parameter choices.
This choice may simultaneously satisfy nonvanishing conditions for
any finite collection of nonzero real polynomials.
\end{lemma}
\begin{proof}
The column spans of different connected components form a direct sum,
so the identity $\sum_i\vect{\alpha}_i\otimes\mathcal D(\vect{t}_i)=0$
holds separately on each component. In a real basis of $\set{B}$, each
coordinate of the subsum is a real vector killed by
$\mathcal D$, proving membership in $\set{B}\otimes \set{V}_0$.
Let $\matr{C}_{\set{I}}$ be the column matrix of $\mathcal D(\vect{t}_i)$ and $\matr{A}_{\set{I}}$
the coefficient matrix of $\vect{\alpha}_i$. The component relation
gives $\matr{C}_{\set{I}}\matr{A}_{\set{I}}^{\mathsf T}=0$, so
$\operatorname{im}\matr{A}_{\set{I}}^{\mathsf T}
\subseteq\ker\matr{C}_{\set{I}}.$ 
Thus
\[
k_{\set{I}}=\operatorname{rank}\matr{A}_{\set{I}}
\le\dim\ker \matr{C}_{\set{I}}=|\set{I}|-c_{\set{I}}.
\]
For the full column matrices $\matr{C}$ and $\matr{A}$, we similarly have $\matr{C}\matr{A}^\top=0$. Since the column spans of the components form a direct sum, $\operatorname{rank}\matr{C}=\sum_{\set{I}} c_{\set{I}}$.
Hence $k=\operatorname{rank}\matr{A}\le N-\sum_{\set{I}} c_{\set{I}}$.

For the polynomial family, consider every sublist of the fixed
decomposition whose subsum does not belong to $\set{B}\otimes\set{V}_0$.
At least one coordinate of its detected image is a nonzero real
polynomial in the parameters.  Take the product of one such
polynomial for each of these finitely many sublists, together with
any additional nonzero polynomials to be retained.
This product is nonzero and hence is nonzero at some real parameter
choice: induction on the number of variables reduces this fact to
the finiteness of the roots of a nonzero univariate polynomial.
At that choice, a subsum with zero detected image belongs to
$\set{B}\otimes\set{V}_0$.  Apply the first part to the selected real map.
In particular, finitely many attainable matrix-rank bounds can be
imposed simultaneously by retaining suitable nonzero minors.
\end{proof}

Connectedness also plays a central role in the splitting theorem for
product tensors of Lovitz and Petrov \cite{lovitz2023generalization}.
The following exterior-support estimate and its equality case follow
from the generic-section formula for subspace arrangements.
For a vector space $\set{W}$, an element of $\bigwedge^d\set{W}$ is
\emph{decomposable} if it has the form
$\vect{w}_1\wedge\cdots\wedge\vect{w}_d$ with $\vect{w}_i\in\set{W}$.
The \emph{support of an exterior vector} is the smallest subspace whose
exterior power contains that vector.

\begin{lemma}[Connected exterior support]
\label{lem:connected-support}
Let $\set{W}$ be a real vector space and let $d\ge1$.
Let $\vect{\omega}_1,\ldots,\vect{\omega}_N\in\bigwedge^d\set{W}$
be nonzero decomposable vectors with connected vector matroid of rank $c$.
Put $\set{L}_i\eqdef\supp\vect{\omega}_i$.
Then
\begin{equation}\label{eq:connected-support}
\dim\sum_i\set{L}_i\le c+d-1.
\end{equation}
If equality holds, every nonempty sublist $\set{I}$ satisfies
$$\dim\sum_{i\in\set{I}}\set{L}_i
\ge\dim\Span_{i\in\set{I}}\vect{\omega}_i+d-1.$$
\end{lemma}

\begin{proof}
Choose generic real covectors $\ell_1,\ldots,\ell_{d-1}$ on $\sum_i\set{L}_i$
and put $\set{U}\eqdef\bigcap_j\ker\ell_j$.
Let $\Phi\eqdef\iota_{\ell_{d-1}}\cdots\iota_{\ell_1}$ denote
successive exterior contraction and put $\vect{v}_i\eqdef\Phi\vect{\omega}_i$.
For $d=1$, take $\set{U}=\sum_i\set{L}_i$ and $\Phi=\operatorname{id}$.
Generically $\RR\vect{v}_i=\set{L}_i\cap\set{U}$. 
Connectedness gives a nonzero vector $\vect{\eta}_{\set{I}}$ in the
intersection of the original column spans across every nonempty proper cut
$\set{I}\mid\set{I}^c$.
The contraction of any nonzero $d$-vector is a nonzero vector-valued
polynomial in the covectors.
Thus we can require $\Phi\vect{\eta}_{\set{I}}\ne0$ for all these
finitely many cuts.
We also choose $\set{U}$ \emph{transverse} to every sublist support sum,
meaning that
$$\dim((\sum_{i\in\set{I}}\set{L}_i)\cap\set{U})
=\dim\sum_{i\in\set{I}}\set{L}_i-d+1$$ 
for every nonempty sublist $\set{I}$.
These conditions hold simultaneously by avoiding the zero sets of
finitely many nonzero polynomials in the covectors.
For every cut, $\Phi\vect{\eta}_{\set{I}}$ is a nonzero vector in
both projected column spans, so these spans cannot form a direct sum.
Thus the vector matroid of the $\vect{v}_i$ is connected.

The generic-section formula
\cite[Theorems~5.5 and~7.1]{raz2019subspace} gives
\[
\dim\Span_i\vect{v}_i
=\min_{\Pi}\sum_{\set{J}\in\Pi}
\left(\dim\sum_{i\in\set{J}}\set{L}_i-d+1\right),
\]
where $\Pi$ ranges over partitions of $\{1,\ldots,N\}$.
For a minimizing partition, transversality gives
$$\dim\Span_i\vect{v}_i\le\sum_{\set{J}\in\Pi}\dim\Span_{i\in\set{J}}\vect{v}_i
\le\sum_{\set{J}\in\Pi}(\dim\sum_{i\in\set{J}}\set{L}_i-d+1)=\dim\Span_i\vect{v}_i.$$
Hence its column spans form a direct sum.
Connectedness forces the one-part partition, so
$$\dim\Span_i\vect{v}_i=\dim\sum_i\set{L}_i-d+1\le c,$$ proving \eqref{eq:connected-support}.
In the equality case, $\dim\Span_i\vect{v}_i=c
=\dim\Span_i\vect{\omega}_i$, so $\Phi$ is injective on
$\Span_i\vect{\omega}_i$.
For every nonempty $\set{I}$, transversality therefore gives
$$\dim\Span_{i\in\set{I}}\vect{\omega}_i
=\dim\Span_{i\in\set{I}}\vect{v}_i
\le\dim((\sum_{i\in\set{I}}\set{L}_i)\cap\set{U})
=\dim\sum_{i\in\set{I}}\set{L}_i-d+1.$$
The proof is complete. 
\end{proof}

\subsection{Hankel decomposition blocks}

We use one Hankel detector for every number $q\ge2$ of detected modes.
Assume $n\ge3$ and take
$\RR^{n-1}\oplus(\RR^n)^{\oplus(q-2)}$ as the ambient space.
For $\vect{x}=(x_1,\ldots,x_n)$, write
$$\vect{x}^-\eqdef(x_1,\ldots,x_{n-1})\ \ \textrm{and}\ \ 
\vect{x}^+\eqdef(x_2,\ldots,x_n).$$ 
Let $\matr{L}_-\vect{x}\eqdef\vect{x}^-$ and
$\matr{L}_+\vect{x}\eqdef\vect{x}^+$.
For real parameters $a_j,b_j,c_{j\ell}$ with
$1\le j\le q$ and $1\le\ell\le q-2$, define
\begin{equation}\label{def:Dmx-hyb}
\begin{aligned}
\matr{A}_j^{\mathrm h}\vect{x}
\eqdef(a_j\vect{x}^-+b_j\vect{x}^+,
c_{j1}\vect{x},\ldots,c_{j,q-2}\vect{x});\ \
\mathcal D_{\mathrm h}(\vect{x}_1\otimes\cdots\otimes\vect{x}_q)
\eqdef\matr{A}_1^{\mathrm h}\vect{x}_1\wedge\cdots\wedge
\matr{A}_q^{\mathrm h}\vect{x}_q.
\end{aligned}
\end{equation}
For $q=2$, only the truncated block occurs.
For $\set{S}\subseteq\RR^n$, put
$\set{S}^{\pm}\eqdef\{\vect{v}^{\pm}:\vect{v}\in\set{S}\}$ and
$\delta_{\mathrm{sh}}(\set{S})\eqdef\dim(\set{S}^-+\set{S}^+)$. 
For $j\ge0$, put $\vect{v}_j(t)\eqdef(1,t,\ldots,t^j)$.
Put $\vect{v}_{n-1}(\infty)\eqdef(0,\ldots,0,1)$ for the direction
of the \emph{node at infinity}.
A \emph{node power} is a nonzero real multiple of
$\vect{v}_{n-1}(\theta)^{\otimes q}$ for
$\theta\in\RR\cup\{\infty\}$.
For the auxiliary case $q=2$, $\set{H}_2(n)$ denotes the space
of real Hankel matrices.

\begin{lemma}[Hankel kernel and exact support expansion]
\label{lem:hybrid-kernel-expansion}
The common kernel of \eqref{def:Dmx-hyb} is $\set{H}_q(n)$.
For every fixed $\set{S}\subseteq\RR^n$, outside the zero set of
a nonzero real polynomial in the parameters, we have
\begin{equation}\label{eq:hybrid-expansion}
\dim\sum_{j=1}^q\matr{A}_j^{\mathrm h}\set{S}
=(q-2)\dim\set{S}+\delta_{\mathrm{sh}}(\set{S}).
\end{equation}
If $\set{S}\ne\RR^n$, then $\delta_{\mathrm{sh}}(\set{S})\ge\dim\set{S}$.
For $\set{S}=\RR^n$, its value is $n-1$.
\end{lemma}

\begin{proof}
All images $\matr{A}_j^{\mathrm h}\vect{v}_{n-1}(t)$ lie in
$\RR\vect{v}_{n-2}(t)\oplus
(\RR\vect{v}_{n-1}(t))^{\oplus(q-2)}$.
This space has dimension $q-1$, so its $q$-fold wedge vanishes.
Node powers at $q(n-1)+1$ distinct real nodes span $\set{H}_q(n)$,
which therefore lies in the common kernel.
Conversely, specialize two chosen modes to the truncated block,
using independently the minus and plus shifts.
Assign the remaining $q-2$ modes to distinct full blocks.
Taking coordinates there shows that every two-mode matrix fiber
$\matr{M}$ satisfies
\[
\matr{L}_-\matr{M}\matr{L}_-^{\mathsf T},\
\matr{L}_-\matr{M}\matr{L}_+^{\mathsf T},\
\matr{L}_+\matr{M}\matr{L}_-^{\mathsf T},\
\matr{L}_+\matr{M}\matr{L}_+^{\mathsf T}
\ \text{are symmetric}.
\]
The first and fourth matrices give symmetry away from the opposite
corners. The mixed equations give $M_{1,n}=M_{n-1,2}$ and
$M_{n,1}=M_{2,n-1}$, whose right-hand entries agree.
Thus $\matr{M}$ is symmetric.
The second matrix now gives $M_{i,j+1}=M_{i+1,j}$.
Every two-mode fiber is therefore Hankel.
Moving one unit between indices connects all multiindices of a fixed
total, so the tensor is Hankel.

The image in \eqref{eq:hybrid-expansion} lies in
$(\set{S}^-+\set{S}^+)\oplus\set{S}^{\oplus(q-2)}$.
Specialize two maps to the minus and plus shifts in the truncated
block and the other maps to identities in distinct full blocks.
Their images span this entire space.
A maximal minor at this specialization gives the required nonzero
polynomial, proving \eqref{eq:hybrid-expansion}.

For proper $\set{S}$, choose a nonzero covector $\ell$ vanishing on it.
The kernel of $\matr{L}_-+z\matr{L}_+$ for $z\ne0$ is spanned by
$\vect{k}(z)\eqdef(z^{n-1},-z^{n-2},\ldots,(-1)^{n-1})$.
The polynomial $\ell(\vect{k}(z))$ is nonzero.
Choose $z\ne0$ at which it does not vanish.
Then $\matr{L}_-+z\matr{L}_+$ is injective on $\set{S}$, giving
$\delta_{\mathrm{sh}}(\set{S})\ge\dim\set{S}$.
The full-space value is immediate.
\end{proof}

\begin{lemma}[Hankel component subsums]
\label{lem:hankel-components}
Let $0\ne\tens{Y}\in\set{B}\otimes(\RR^n)^{\otimes q}$ be Hankel
in its $q\ge2$ tail modes, where $n\ge3$.
Fix a shortest CP decomposition
\begin{equation*}
\tens{Y}=\sum_{i=1}^N\vect{\alpha}_i\otimes
\bigotimes_{j=1}^q\vect{x}_{ij},
\end{equation*}
and put $k\eqdef\dim\Span_i\vect{\alpha}_i$.
Choose parameters as in \cref{lem:detector-components}, imposing
\eqref{eq:hybrid-expansion} for every mode-factor span of every sublist.
For each connected component $\set{I}$ of nonzero detector columns,
its subsum is tail-Hankel of rank $|\set{I}|$ with common tail support
$\set{S}_{\set{I}}$.
Write $c_{\set{I}}$ for its column rank and
$k_{\set{I}}\eqdef\dim\Span_{i\in\set{I}}\vect{\alpha}_i$.
Then
\begin{equation}\label{eq:component-budget}
k_{\set{I}}+c_{\set{I}}\le|\set{I}|;
\ \ \sum_{\set{I}}c_{\set{I}}\le N-k,
\end{equation}
and
\begin{equation}\label{eq:detector-component-support}
(q-2)\dim\set{S}_{\set{I}}+\delta_{\mathrm{sh}}(\set{S}_{\set{I}})
\le c_{\set{I}}+q-1\le|\set{I}|-k_{\set{I}}+q-1.
\end{equation}
Zero-column singletons have node-power tails, including the node at infinity.
\end{lemma}

\begin{proof}
By \cref{lem:hybrid-kernel-expansion}, the detector annihilates
$\tens{Y}$ for every parameter choice.
Thus \cref{lem:detector-components} puts each component subsum in
the common tail kernel and gives \eqref{eq:component-budget}.
By \cref{lem:shortest-subsum-support}, the subsum has rank $|\set{I}|$
and its displayed factor spans are its intrinsic supports.
Tail symmetry identifies these with $\set{S}_{\set{I}}$.
For nonzero columns, their exterior supports span
$\sum_j\matr{A}_j^{\mathrm h}\set{S}_{\set{I}}$.
Combining \eqref{eq:hybrid-expansion} with
\cref{lem:connected-support} proves \eqref{eq:detector-component-support}.
A zero column gives a simple tail in the common kernel.
It is symmetric, and its rank-one Hankel matrix identities force
its direction to be proportional to $(1,t,\ldots,t^{n-1})$
or $(0,\ldots,0,1)$.
\end{proof}

Apply the detector \eqref{def:Dmx-hyb} with $q=m$ and scalar
coefficient space. The component estimate then takes the following form.

\begin{proposition}[Hybrid component bound]
\label{prop:hybrid-component}
Let $0\ne\tens{H}\in\set{H}_m(n)$ with $m\ge3$ and $n\ge3$, and let
$\tens{H}=\sum_{i=1}^{N}\tens{t}_i$ be a shortest real CP decomposition.
With parameters chosen as in \cref{lem:detector-components} so that
\eqref{eq:hybrid-expansion} holds for every mode-factor span of every sublist,
let $\set{I}$ be a connected component of the nonzero columns $\mathcal D_{\mathrm h}(\tens{t}_i)$, with column rank
$c_{\set{I}}$.  Let $\set{S}_{\set{I}}$ be the common mode support of its
subsum and put $e_{\set{I}}=\dim\set{S}_{\set{I}}$.  Then
\begin{equation}\label{eq:hybrid-component-bound}
(m-2)e_{\set{I}}+\delta_{\mathrm{sh}}(\set{S}_{\set{I}})\le c_{\set{I}}+m-1\le |\set{I}|+m-2.
\end{equation}
Consequently
\begin{equation}\label{eq:hybrid-component-alternatives}
|\set{I}|\ge
\begin{cases}
(m-1)(e_{\set{I}}-1)+1,& e_{\set{I}}<n,\\
(m-1)(n-1),& e_{\set{I}}=n.
\end{cases}
\end{equation}
Every term corresponding to a zero column is a nonzero real multiple of a node power, including the node at infinity.
\end{proposition}

\begin{proof}
Apply \cref{lem:hankel-components} with $q=m$ and scalar
coefficient space, so $k_{\set{I}}=1$.
Equation \eqref{eq:detector-component-support} becomes
\eqref{eq:hybrid-component-bound}, and
\cref{lem:hybrid-kernel-expansion} gives
\eqref{eq:hybrid-component-alternatives}.
The zero-column assertion is also part of \cref{lem:hankel-components}.
\end{proof}

\subsection{Operator-balanced and symmetric blocks}

Let $q\ge3$ and let $\set{E}$ be a real vector space.
For a basis $\vect{f}_1,\ldots,\vect{f}_{q-1}$ of $\RR^{q-1}$, define
\[
\matr{B}_j\vect{x}\eqdef\vect{x}\otimes\vect{f}_j
\ (1\le j<q);
\ \matr{B}_q\vect{x}\eqdef
\vect{x}\otimes(\vect{f}_1+\cdots+\vect{f}_{q-1}).
\]
Put $\mathcal D_{\mathrm{sym}}(\bigotimes_{j=1}^q\vect{x}_j)
\eqdef\bigwedge_{j=1}^q\matr{B}_j\vect{x}_j$.
For $\matr{C}\in\operatorname{End}(\set{E})$, let $\matr{C}_j$
act in mode $j$ and set
$$\set{K}_{\matr{C}}\eqdef
\{\tens{T}\in\Sym^q(\set{E}):\matr{C}_1\tens{T}=\cdots=\matr{C}_q\tens{T}\}.$$ 
We call tensors in $\set{K}_{\matr{C}}$ \emph{operator-balanced} with respect to $\matr{C}$.
We will also use the detector
$$\mathcal D_z\eqdef\mathcal D_{\mathrm{sym}}\circ
(\operatorname{id}^{\otimes(q-1)}\otimes(\operatorname{id}+z\matr{C})).$$

\begin{lemma}[Operator-balanced decomposition blocks]
\label{lem:operator-balanced-components}
Let $\tens{Y}\in\set{B}\otimes\set{K}_{\matr{C}}$.
Every shortest CP decomposition
$\tens{Y}=\sum_i\vect{\alpha}_i\otimes\bigotimes_{j=1}^q\vect{x}_{ij}$
partitions, for a suitable real $z$, into detector components
$\tens{Y}_a\in\set{B}\otimes\set{K}_{\matr{C}}$.
Each component has a common tail support $\set{S}_a$ with
$\matr{C}\set{S}_a\subseteq\set{S}_a$.
Its length $N_a=\rk(\tens{Y}_a)$ satisfies
$$N_a\ge(q-1)(\dim\set{S}_a-1)+k_a,$$
where $k_a\eqdef\dim\Span_{i\text{ in component }a}\vect{\alpha}_i$.
Zero columns are singletons whose tails are multiples of
$\vect{v}^{\otimes q}$ with $\matr{C}\vect{v}\in\RR\vect{v}$.
For a shortest signed power decomposition of
$\tens{Y}\in\set{K}_{\matr{C}}$, the same component and support
conclusions hold, and the component lengths $N_a=s(\tens{Y}_a)$ satisfy
$$N_a\ge(q-1)(\dim\set{S}_a-1)+1.$$
\end{lemma}

\begin{proof}
The exterior coordinates with degree two in auxiliary block $j$
and degree one in every other block give antisymmetrization
in modes $j$ and $q$. Conversely, powers have zero image and span
$\Sym^q(\set{E})$. Thus $\ker\mathcal D_{\mathrm{sym}}=\Sym^q(\set{E})$.
The constant and linear coefficients of $\mathcal D_z$ then show
$\bigcap_z\ker\mathcal D_z=\set{K}_{\matr{C}}$.
Indeed, the linear coefficient requires $\matr{C}_q\tens{T}$ to be
symmetric, which for symmetric $\tens{T}$ is equivalent to
$\matr{C}_j\tens{T}=\matr{C}_q\tens{T}$ for every $j$.

Apply \cref{lem:detector-components,lem:shortest-subsum-support}.
The selected components lie in the common kernel, and their
displayed tail-factor spans equal their intrinsic supports.
Tail symmetry identifies these supports with $\set{S}_a$.
Moving $\matr{C}$ to another tail mode before contraction gives
$\matr{C}\set{S}_a\subseteq\set{S}_a$.
For a nonzero-column component of column rank $c_a$, the first
$q-1$ exterior factors span $\set{S}_a\otimes\RR^{q-1}$.
The last factors lie in this space by invariance.
Hence \cref{lem:connected-support,lem:detector-components} give
$$(q-1)\dim\set{S}_a\le c_a+q-1\le N_a-k_a+q-1.$$
A nonzero simple tensor in the common kernel is a multiple of
$\vect{v}^{\otimes q}$ with $\matr{C}\vect{v}\in\RR\vect{v}$.
Its singleton satisfies $N_a=\dim\set{S}_a=k_a=1$.
For signed powers, use the symmetric part of
\cref{lem:shortest-subsum-support}, with $N_a=s(\tens{Y}_a)$ and $k_a=1$;
the same argument gives the asserted length bound.
\end{proof}

\begin{lemma}[Symmetric decomposition blocks]
\label{lem:symmetric-blocks}
Every shortest CP decomposition of $\tens{S}\in\Sym^m(\RR^n)$
with $m\ge3$ and $n\ge1$
partitions into symmetric subsums $\tens{F}_a\in\Sym^m(\set{E}_a)$
of lengths $N_a$, with common intrinsic mode supports
$\set{E}_a\subseteq\RR^n$ satisfying
\begin{equation}\label{eq:block-density}
\rk(\tens{F}_a)=N_a\ge(m-1)(\dim\set{E}_a-1)+1.
\end{equation}
If $\tens{S}$ is concise and $\rk(\tens{S})\le(m-1)(n-1)$,
every $\set{E}_a$ is proper.
\end{lemma}

\begin{proof}
Apply \cref{lem:operator-balanced-components} with
$q=m$ and $\set{E}=\RR^n$, using $\set{B}=\RR$ and $\matr{C}=0$.
Then $k_a=1$, giving \eqref{eq:block-density}.
A full-support component would have length at least
$$(m-1)(n-1)+1>\rk(\tens{S}),$$ proving the last assertion.
\end{proof}

\section{Restriction estimates and CP rank additivity}
\label{sec:restriction}


\subsection{Coprime recurrence groups and polynomial projections}
\label{subsec:window-recurrence}

We merge components with a common recurrence factor and estimate their
lengths in the original polynomial window.

\begin{lemma}[Coprime recurrence groups]
\label{lem:window-recurrence-groups}
Let $m\ge3$ and $n\ge3$. Let $\set{K}$ be a real vector space and let
$0\ne\tens{Y}\in\set{K}^*\otimes(\set{P}_{n-1}^*)^{\otimes(m-1)}$
be Hankel in its tail modes.
Suppose the components of a shortest decomposition supplied by
\cref{lem:hankel-components} all have proper common tail supports.
Let $f_a$ be their recurrences from \cref{lem:family-recurrence}.
Choose a common real affine chart preserving their degrees and use
monic dehomogenizations. Put $f\eqdef\lcm_a f_a$ and assume $\deg f\le M$.
Let $p\mid f$ be monic and nonconstant and let $\lambda$ be Frobenius on
$\set{A}_p\eqdef\RR[t]/(p)$.
Suppose a linear map $\iota:\set{K}\to\set{A}_p$ satisfies
\[
\tens{Y}(v,w_1,\ldots,w_{m-1})
=\lambda\left(\iota(v)\prod_{j=1}^{m-1}[w_j]_p\right)
\ (v\in\set{K};\ w_j\in\set{P}_{n-1}).
\]
Then $p=\prod_{i=1}^b p_i$ for pairwise coprime nonconstant polynomials
$p_i$ such that
\begin{equation}\label{eq:window-recurrence-groups}
\rk(\tens{Y})\ge(m-2)\deg p+\sum_{i=1}^b(k_i-1)-(m-3)b,
\end{equation}
where $k_i\eqdef\dim\pi_i(\iota(\set{K}))$ for the quotient
$\pi_i:\set{A}_p\to\set{A}_{p_i}$.
In particular $b\le\tau(p)$ and $k_i\ge1$.
\end{lemma}

\begin{proof}
For an initial component $\tens{Y}_a$ of length $N_a$ and common
tail support $\set{S}_a$, we have $\deg f_a=\dim\set{S}_a$.
Put $k_a\eqdef\dim\supp_1(\tens{Y}_a)$.
By \cref{lem:hybrid-kernel-expansion,lem:hankel-components},
we have $$N_a\ge(m-2)(\deg f_a-1)+k_a.$$
This also holds for zero-column singletons.
Merge two components whenever their recurrences have a common factor,
using their least common multiple as the new recurrence.
Its degree is at most the sum of the two degrees minus one,
and the first-support dimension is at most the sum of the two
first-support dimensions. Thus the preceding length bound is
preserved at every merge.
At termination the recurrences $q_i$ are pairwise coprime with
$\prod_iq_i=f$. If $N_i$ is the merged length and $\kappa_i$ its
first-support dimension, then
\[
N_i\ge(m-2)(\deg q_i-1)+\kappa_i.
\]

Division by $q_i$ expresses the moment functional of the $i$th
merged subsum as a functional
$\mu_{i,v}\in\set{A}_{q_i}^*$.
The products of $m-1$ window polynomials span $\RR[t]_{\le M-1}$.
Since $\deg f\le M$, their classes span $\set{A}_f$.
Consequently the given identity for $\tens{Y}$ determines the sum
of the $\mu_{i,v}$ on this whole quotient.
The Chinese remainder decomposition makes its components unique.
Specifically, for the quotient $\pi:\set{A}_f\to\set{A}_p$,
the sum is the pullback $u\mapsto\lambda(\iota(v)\pi(u))$.

Put $p_i\eqdef\gcd(p,q_i)$.
Then $\prod_i p_i=p$, and uniqueness identifies the $i$th component
with the functional induced by $\lambda$ on the $p_i$-factor.
A constant $p_i$ would make the merged subsum zero, contradicting
\cref{lem:shortest-subsum-support}.
For each nonconstant $p_i$, let $\lambda_i$ be the Frobenius
functional on $\set{A}_{p_i}$ induced by the Chinese remainder
decomposition of $\set{A}_p$.
The merged moment functional is therefore the pullback of
$z\mapsto\lambda_i(\pi_i(\iota(v))z)$.
Because $\deg p_i\le M$, window products cover $\set{A}_{p_i}$.
Nondegeneracy of $\lambda_i$ now gives
$\kappa_i=\dim\pi_i(\iota(\set{K}))=k_i$.
Summing the length bounds and using $\deg p_i\le\deg q_i$
proves \eqref{eq:window-recurrence-groups}.
Coprimality gives $b\le\tau(p)$, and the nonzero merged subsums give $k_i\ge1$.
\end{proof}

The following block-count estimate is the real version of
\cite[Lemma~6.1]{li2026hankeladditivity}; its proof works over $\RR$
without change.

\begin{lemma}[Block count for a polynomial subspace]
\label{lem:polynomial-block-count}
Let $p=\prod_{i=1}^b p_i\in\RR[t]$ be monic of degree $D\ge1$ with pairwise
coprime nonconstant factors.
Let $c\in\mathbb Z_{\ge0}$ satisfy $c<D$ and let $0\ne\set{U}\subseteq\RR[t]_{\le D-c-1}$.
Let $\pi_i:\set{A}_p\to\set{A}_{p_i}$ denote reduction.
Then we have 
$$\sum_i\dim\pi_i(\set{U})\ge c+\dim\set{U}\ \ \textrm{or}\ \ b\le D-\dim\set{U}.$$ 
\end{lemma}

\begin{proof}
Suppose $\sum_i\dim\pi_i(\set{U})<c+\dim\set{U}$.
On $\set{A}_p\eqdef\RR[t]/(p)$, the coefficient functional
$\chi([z]_p)\eqdef[t^{D-1}][z]_p$ is Frobenius.
Identify polynomial subspaces with their residue classes and put
$\set{B}\eqdef\bigoplus_i\pi_i(\set{U})$ under the Chinese remainder
isomorphism. Write
$\set{B}^{\perp_\chi}\eqdef\{a\in\set{A}_p:\chi(ab)=0\text{ for every }b\in\set{B}\}$.
Since $\chi(\RR[t]_{\le c-1}\set{U})=0$, requiring orthogonality
to $\set{B}$ imposes at most $\dim\set{B}-\dim\set{U}$
additional linear conditions on $\RR[t]_{\le c-1}$.
Thus the space
$\set{Z}\eqdef\RR[t]_{\le c-1}\cap\set{B}^{\perp_\chi}$
satisfies 
$$\dim\set{Z}\ge c+\dim\set{U}-\dim\set{B}>0.$$
Put $q\eqdef\gcd(p,\set{Z})$.
The inclusion $\set{Z}\subseteq q\RR[t]_{\le c-1-\deg q}$
gives $\dim\set{Z}\le c-\deg q$, hence 
$$\deg q\le c-\dim\set{Z}.$$ 
If $\pi_i(\set{Z})=0$, then $p_i\mid q$.
Thus at most $c-\dim\set{Z}$ projections vanish.
The Chinese remainder factors are orthogonal for $\chi$, so
$\bigoplus_i\pi_i(\set{Z})\subseteq\set{B}^{\perp_\chi}$.
Nondegeneracy gives
\[
D-\dim\set{B}
\ge\sum_i\dim\pi_i(\set{Z})
\ge b-c+\dim\set{Z}\ge b+\dim\set{U}-\dim\set{B}.
\]
Hence $b\le D-\dim\set{U}$.
\end{proof}

\subsection{A common symmetric rank restriction estimate}
\label{subsec:symmetric-restriction}

\begin{proposition}[Symmetric rank under a short restriction]
\label{prop:symmetric-restriction}
Let $\tens{H}\in\set{H}_m(n)$ be concise with $m\ge3$ and $n\ge3$.
Let $0\ne\set{K}\subsetneq\set{P}_{n-1}$ and
put $h\eqdef\gcd\set{K}$.
Suppose $$\rk(\tens{H}_{\set{K}})<(m-1)\dim\set{K}\ \ \textrm{and}\ \ 
\rk(\tens{H}_{\set{K}})-\dim\set{K}<(m-2)(n-1)-1.$$ 
Then
\begin{equation*}
s(\tens{H})\le\rk(\tens{H}_{\set{K}})+\codim\set{K}+(m-2)\deg h.
\end{equation*}
\end{proposition}

\begin{proof}
Concision of $\tens{H}$ makes $\tens{H}_{\set{K}}$ concise in its first input.
Apply \cref{lem:hankel-components} with $q=m-1$ and first-support
dimension $\dim\set{K}$.
By \eqref{eq:component-budget}--\eqref{eq:detector-component-support},
a full-tail-support component would require
$$(m-2)(n-1)-1\le c_a\le\rk(\tens{H}_{\set{K}})-\dim\set{K},$$ contrary to the hypothesis.
Thus every component has proper tail support.
By \cref{lem:family-recurrence}, its recurrence $f_a$ has degree
equal to that support dimension and hence at most its length.
The form $f\eqdef\lcm_a f_a$ therefore satisfies $\deg f\le\rk(\tens{H}_{\set{K}})$ and
$L(\set{K}f\set{P}_{(m-1)(n-1)-\deg f})=0$, where $L$ is the moment
functional of $\tens{H}$.
Since $\deg h\le\codim\set{K}$ and
$$\deg f+\codim\set{K}-\deg h\le(m-1)\dim\set{K}-1+\codim\set{K}-\deg h\le(m-1)(n-1),$$
\cref{lem:apolar-lifting} makes $hf$ apolar.
Choose common real binary coordinates preserving $\deg h$ and
all $\deg f_a$ and use monic dehomogenizations.
By \cref{lem:apolar-quotient}, a divisor $g\mid hf$ and a Frobenius
functional $\lambda_g$ on $\set{A}_g\eqdef\RR[t]/(g)$ represent
$\tens{H}$ with $s(\tens{H})\le\Gamma_m(g)$.
Concision forces $\deg g\ge n$.

Put $h_0\eqdef\gcd(g,h)$.
Then $g/h_0$ is nonconstant and $(g/h_0)\mid f$, since
$$\deg h<n\le\deg g\ \ \textrm{and}\ \  \gcd(g/h_0,h/h_0)=1.$$
The functional
$\lambda_0([z]_{g/h_0})\eqdef\lambda_g([h_0z]_g)$ is well-defined.
If $\lambda_0([z]_{g/h_0}\set{A}_{g/h_0})=0$, the quotient map
$\set{A}_g\to\set{A}_{g/h_0}$ and nondegeneracy of $\lambda_g$
give $g\mid h_0z$, hence $(g/h_0)\mid z$.
Thus $\lambda_0$ is Frobenius.
Then $\dim(\set{K}/h_0)=\dim\set{K}$ and
$$\set{K}/h_0\subseteq\RR[t]_{\le n-\deg h_0-1}.$$
Reduction on $\set{K}/h_0$ is injective:
if $(g/h_0)\mid v/h_0$ for $v\in\set{K}$, then $g\mid v$;
the contraction of $\tens{H}$ by $v$ vanishes, so concision gives $v=0$.
The window identity is
\[
\tens{H}_{\set{K}}(v,w_1,\ldots,w_{m-1})
=\lambda_0\left([v/h_0]_{g/h_0}\prod_{j=1}^{m-1}[w_j]_{g/h_0}\right).
\]

Apply \cref{lem:window-recurrence-groups} with
$\iota(v)\eqdef[v/h_0]_{g/h_0}$, using $\deg f\le\rk(\tens{H}_{\set{K}})\le(m-1)(n-1)$.
For its factors $p_i$, we have
$\sum_i k_i-\dim\set{K}\ge0$.
Its bound becomes
\[
\rk(\tens{H}_{\set{K}})
\ge(m-2)(\deg g-\deg h_0-b)+\sum_i k_i.
\]
Apply \cref{lem:polynomial-block-count} with polynomial $g/h_0$,
subspace $\set{K}/h_0$, and parameter $\deg g-n$.
The polynomial has degree $\deg g-\deg h_0$, whose difference from
the parameter is $n-\deg h_0>0$, giving the required window degree.
Thus either $$\sum_i k_i-\dim\set{K}\ge\deg g-n\ \ 
\textrm{or}\ \
b\le\deg g-n+\codim\set{K}-\deg h_0.$$ 
The latter alternative would give $\rk(\tens{H}_{\set{K}})\ge(m-1)\dim\set{K}$, contrary to the
hypothesis. Hence $\sum_i k_i-\dim\set{K}\ge\deg g-n$.
Substituting this into the preceding bound and using
$b\le\tau(g/h_0)\le\tau(g)$ gives
\begin{align*}
\rk(\tens{H}_{\set{K}})
&\ge(m-2)(\deg g-\deg h_0-b)+\dim\set{K}+\deg g-n,\\
s(\tens{H})&\le(m-1)\deg g-(m-2)\tau(g)\\
&\le\rk(\tens{H}_{\set{K}})+\codim\set{K}+(m-2)\deg h_0
\le\rk(\tens{H}_{\set{K}})+\codim\set{K}+(m-2)\deg h.
\end{align*}
The last inequality uses $h_0\mid h$.
\end{proof}

\subsection{CP rank additivity}



Let $m\ge3$ and let
$\tens{T}\in\set{A}_1\otimes\cdots\otimes\set{A}_m$.
For $1\le i\le m$ and a linear subspace
$\set{C}\subseteq\set{A}_i$, denote by
$\tens{T}_{i,\set{C}}$ the image of $\tens{T}$ under the quotient
$\set{A}_i\to\set{A}_i/\set{C}$ in mode $i$.
We say that $\tens{T}$ satisfies the
\emph{quotient rank-loss condition} if there exists a set
$\set{I}\subseteq\{1,\ldots,m\}$ with $|\set{I}|=m-1$
such that, for every $i\in\set{I}$ and every linear subspace
$\set{C}\subseteq\set{A}_i$, we always have 
\begin{equation}\label{eq:quotient-rank-loss}
\rk(\tens{T})-\rk(\tens{T}_{i,\set{C}})
\le(m-1)\dim\set{C}.
\end{equation}
The following is the real version of
\cite[Proposition~3.3]{li2026hankeladditivity}.
Its proof applies without change over $\RR$, since it uses only
linear quotients, the matroid union theorem, and extensions of
linear maps.

\begin{lemma}[Restriction and elimination criterion]
\label{prop:restriction-criterion}
Let $m\ge3$ and let
$\tens{T}\in\set{A}_1\otimes\cdots\otimes\set{A}_m$
satisfy the quotient rank-loss condition in \eqref{eq:quotient-rank-loss}.
Then $\Delta_r(\tens{T},\tens{X})=0$ for every real tensor
$\tens{X}$ of order $m$.
\end{lemma}

The preceding component and recurrence estimates now give the
restriction bound needed to verify the quotient rank-loss condition in \eqref{eq:quotient-rank-loss}.  

\begin{proposition}[Restriction inequality for real Hankel tensors]
\label{thm:restriction-hankel}
Let $\tens{H}\in\set{H}_m(n)$ be concise with $m\ge3$ and $n\ge3$,
and let $L\in\set{P}_{m(n-1)}^*$ be its moment functional.
For every linear subspace $\set{K}\subseteq \set{P}_{n-1}$ and
$\tens{H}_{\set{K}}\eqdef\tens{H}|_{\set{K}\times \set{P}_{n-1}^{m-1}}$, we have
\begin{equation*}
\rk(\tens{H})-\rk(\tens{H}_{\set{K}}) \le (m-1)\codim \set{K}.
\end{equation*}
\end{proposition}

\begin{proof}
The cases $\set{K}=0$ and $\set{K}=\set{P}_{n-1}$ are immediate.
The universal bound \cref{cor:ordinary-upper} also settles the claim
when $\rk(\tens{H}_{\set{K}})+(m-1)\codim\set{K}\ge M$.
Otherwise $$\rk(\tens{H}_{\set{K}})\le(m-1)(\dim\set{K}-1).$$ 
Since $1\le\dim\set{K}\le n-1$, both hypotheses of
\cref{prop:symmetric-restriction} hold.
Using $\deg\gcd\set{K}\le \codim\set{K}$, we obtain
$$\rk(\tens{H})\le s(\tens{H})\le \rk(\tens{H}_{\set{K}})+\codim\set{K}+(m-2)\codim\set{K}=\rk(\tens{H}_{\set{K}})+(m-1)\codim\set{K}.$$ 
The proof is complete. 
\end{proof}

\begin{proof}[Proof of \cref{thm:main-additivity}]
The rank bound is \cref{cor:ordinary-upper}.
The case $\tens{H}=0$ is immediate.
For a nonzero nonconcise $\tens{H}$,
\cref{prop:nonconcise-hankel} gives a complete Frobenius moment tensor
$\tens{F}$ and a surjective window map.
The tensor $\tens{F}$ is concise and Hankel in the monomial basis.
The window map and a linear section, extended by the identity on
the partner spaces, give
$$\rk(\tens{H}\oplus\tens{T})=\rk(\tens{F}\oplus\tens{T})\ \ \textrm{and}\ \  \rk(\tens{H})=\rk(\tens{F}).$$ 
Thus it suffices to treat concise Hankel tensors.

For $n\ge3$, symmetry and \cref{thm:restriction-hankel} verify
the quotient rank-loss condition through the
quotient--restriction correspondence.
For $n=2$, every line restriction has rank at least one and
$\rk(\tens{H})\le m$; the zero restriction is covered by
$m\le2(m-1)$.
For $n=1$, the rank is one and the condition is immediate.
Now \cref{prop:restriction-criterion} proves \eqref{eq:main-additivity}.
\end{proof}

\section{Rank equalities for real Hankel tensors}
\label{sec:critical}

\subsection{Low-rank equality}
\label{sec:short-rank}

\begin{proposition}[Uniform short rank equality]
\label{thm:short-comon}
Let $\tens{H}\in\set{H}_m(n)$ with $m\ge3$ and $n\ge1$, and let $M$ be as in
\eqref{eq:M_value}.  Then
\[
\rk(\tens{H})\le M-2
\ \Longrightarrow\ 
\rk(\tens{H})=p_2(\tens{H})
=p_{m-1}(\tens{H})=s(\tens{H}).
\]
\end{proposition}

\begin{proof}
The cases $\tens{H}=0$ and $n=1$ are immediate.
For $n=2$, apply
\cite[Corollary~3.12]{zhang2016comon} over $\RR$. 
Assume $n\ge3$ and $\rk(\tens{H})\le M-2$.  By
\eqref{eq:hybrid-component-alternatives}, a full-support component of
a shortest decomposition would require at least $M-1$ terms.  Therefore every
nonzero component has proper support.  It is a nonconcise Hankel
tensor and, by \cref{prop:nonconcise-hankel}, can be replaced by
exactly the same number of signed real powers.  
Terms corresponding to zero columns are already signed real powers. 
Concatenating all replacements produces a
real symmetric decomposition of length $\rk(\tens{H})$.  The rank chain
\eqref{eq:rank-chain} forces equality of all ranks under consideration.
\end{proof}

\subsection{Critical circuits and weak density}
\label{subsec:critical-density}

\begin{lemma}[Structure of a possible endpoint exception]
\label{cor:critical-circuit}
Let $\tens{H}\in\set{H}_m(n)$ with $m\ge3$ and $n\ge3$, and let $M$
be as in \eqref{eq:M_value}. Suppose $\rk(\tens{H})=M-1<s(\tens{H})$.
In every shortest real CP decomposition
$\tens{H}=\sum_{i=1}^{M-1}\tens{t}_i$, choose hybrid parameters
as in \cref{prop:hybrid-component}.
All columns are nonzero and form one \emph{circuit} of rank $M-2$;
that is, they are linearly dependent and every proper sublist is independent.
Their subsum has common mode support $\RR^n$, and
\[
\dim\sum_i\supp\mathcal D_{\mathrm h}(\tens{t}_i)=(m-1)n-1.
\]
\end{lemma}

\begin{proof}
If all components had proper support, the replacement argument in
\cref{thm:short-comon} would give $s(\tens{H})=\rk(\tens{H})$.
Thus a full-support component exists and contains all $M-1$ terms
by \eqref{eq:hybrid-component-alternatives}.
The bounds in \eqref{eq:hybrid-component-bound} force column rank $M-2$.
The relation $\sum_i\mathcal D_{\mathrm h}(\tens{t}_i)=0$
has coefficient vector $(1,\ldots,1)$, which spans the
one-dimensional relation space, so the columns form a circuit.
The support dimension follows from \eqref{eq:hybrid-expansion}.
\end{proof}

\begin{proposition}[Weak density at the critical rank]
\label{lem:hereditary-density}
Let $\tens{H}\in\set{H}_m(n)$ with $m\ge3$ and $n\ge3$, and let $M$
be as in \eqref{eq:M_value}. Suppose $\rk(\tens{H})=M-1<s(\tens{H})$.
Suppose a nonempty proper sublist of a shortest CP decomposition
has a symmetric subsum with support $\set{S}$ and length
\begin{equation}\label{eq:core-length}
N_0\ge(m-1)(\dim\set{S}-1)+1.
\end{equation}
Then $\set{S}\subsetneq\RR^n$ and, for $h\eqdef\gcd(\set{S}^\perp)$,
\begin{equation}\label{eq:gcd-shift-budget}
\dim\set{S}-\deg h\ge1.
\end{equation}
\end{proposition}

\begin{proof}
Properness follows from $N_0<M-1$ and \eqref{eq:core-length}.
Choose hybrid columns $\vect{\omega}_i\eqdef\mathcal D_{\mathrm h}(\tens{t}_i)$
as in \cref{cor:critical-circuit} and put
$\set{L}_i\eqdef\supp\vect{\omega}_i$.
Their rank is $M-2$ and their total support has dimension $M+m-3$.
For the selected proper sublist $\set{I}$, the circuit property
and the equality case of \cref{lem:connected-support} give
$$\dim\sum_{i\in\set{I}}\set{L}_i\ge N_0+m-1.$$ 

We have $\deg h\le\dim\set{S}$.
If equality held, then $\set{S}^\perp=h\set{P}_{n-1-\deg h}$.
The joint annihilator of $\set{S}^-$ and $\set{S}^+$ would be
$$\{p\in\set{P}_{n-2}:h\mid xp,\ h\mid yp\}
=h\set{P}_{n-2-\deg h},$$ giving
$\delta_{\mathrm{sh}}(\set{S})=\dim\set{S}$.
By \cref{lem:shortest-subsum-support} and \eqref{eq:hybrid-expansion},
the selected exterior supports would therefore satisfy
\[
N_0+m-1\le\dim\sum_{i\in\set{I}}\set{L}_i
=(m-1)\dim\set{S},
\]
contrary to \eqref{eq:core-length}.
Thus $\deg h<\dim\set{S}$, proving \eqref{eq:gcd-shift-budget}.
\end{proof}

\subsection{Critical equality and completion}
\label{subsec:higher-order-critical}
\label{subsec:comon-completion}

We apply the strengthened restriction estimate to a symmetric block
satisfying the weak density bound.

\begin{proposition}[Critical equality in every order]
\label{thm:critical-high-order}
Let $\tens{H}\in\set{H}_m(n)$ with $m\ge3$ and $n\ge1$, and let $M$
be as in \eqref{eq:M_value}.
If $\rk(\tens{H})=M-1$, then $s(\tens{H})=\rk(\tens{H})$.
\end{proposition}

\begin{proof}
The cases $\tens{H}=0$ and $n=1$ are immediate.
Nonzero nonconcise tensors are covered by \cref{prop:nonconcise-hankel},
and \cite[Corollary~3.12]{zhang2016comon} over $\RR$ settles $n=2$.
Suppose therefore that $\tens{H}$ is concise with $n\ge3$ and
$\rk(\tens{H})=M-1<s(\tens{H})$.
Partition a shortest decomposition by \cref{lem:symmetric-blocks}.
If every block has dimension at most two, binary rank equality
\cite[Corollary~3.12]{zhang2016comon}
replaces each block by equally many signed powers, a contradiction.
Thus some block has support $\set{S}$ of dimension $e\ge3$
and length $N_0\ge(m-1)(e-1)+1$.
Its support is proper by \cref{lem:symmetric-blocks},
and concision makes it a proper subdecomposition.
Put $\set{K}\eqdef\set{S}^\perp$.
By \cref{lem:hereditary-density}, we have $\deg\gcd\set{K}<e$.
Restriction to $\set{K}$ kills this block, giving
$$\rk(\tens{H}_{\set{K}})\le M-1-N_0\le(m-1)(n-e)-1.$$
Since $e\ge3$, this also gives
$$\rk(\tens{H}_{\set{K}})-\dim\set{K}<(m-2)(n-1)-1.$$
Thus \cref{prop:symmetric-restriction} applies.
Using $\deg\gcd\set{K}\le e-1$, we obtain
\begin{align*}
s(\tens{H})
\le\rk(\tens{H}_{\set{K}})+e+(m-2)\deg\gcd\set{K}
\le M-1-N_0+(m-1)(e-1)+1\le M-1,
\end{align*}
a contradiction.
\end{proof}

\begin{proof}[Proof of \cref{thm:main-comon}]
If $\rk(\tens{H})\le M-1$, then
\cref{thm:short-comon,thm:critical-high-order} give
$s(\tens{H})=\rk(\tens{H})$.
The rank chain \eqref{eq:rank-chain} yields equality with every
partially symmetric rank.
By \cref{cor:ordinary-upper}, only $\rk(\tens{H})=M$ remains.
Then \cref{prop:partial-upper} and the same rank chain give
$$M=\rk(\tens{H})\le p_j(\tens{H})\le p_{m-1}(\tens{H})\le M$$ 
for $1\le j\le m-1$.
\end{proof}

\begin{proof}[Proof of \cref{cor:symmetric-partner}]
By \cref{thm:main-comon} and the hypothesis,
$s(\tens{H})=\rk(\tens{H})$ and $s(\tens{S})=\rk(\tens{S})$.
Concatenating shortest signed power decompositions gives
$s(\tens{H}\oplus\tens{S})\le s(\tens{H})+s(\tens{S})$,
whereas \cref{thm:main-additivity} gives
$$\rk(\tens{H}\oplus\tens{S})=\rk(\tens{H})+\rk(\tens{S}).$$ 
The inequality $\rk\le s$ forces all claimed equalities.
\end{proof}

\section{A fourth-order counterexample at the CP rank bound}
\label{sec:quartic-counterexample}

This section proves the single-block part of \cref{thm:main-counterexample}.
We establish a symmetric rank lower bound and an eleven-power
decomposition attaining it. The CP and partially symmetric ranks then
follow from \cref{prop:partial-upper,thm:main-comon}. 
Let the tensor $\tens{H}_*\in\set{H}_4(4)$ be defined as in  
\eqref{eq:main-counterexample-tensor}. 
Equivalently
\begin{equation*}
\tens{H}_*(f_1,f_2,f_3,f_4)=[t^6](f_1f_2f_3f_4)
\end{equation*}
for $f_i\in\RR[t]_{\le 3}$.
Its associated quartic is
\begin{equation*}
6x_1^2x_4^2+24x_1x_2x_3x_4+4x_1x_3^3
+4x_2^3x_4+6x_2^2x_3^2.
\end{equation*}

\subsection{Quadratic defects and the symmetric rank lower bound}

Throughout this subsection, let $\tens{H}\in\set{H}_4(4)$.
For $\vect{v}=(v_1,v_2,v_3,v_4)\in \RR^4$, use the shifts
$$\vect{v}^-=(v_1,v_2,v_3)\ \ \textrm{and}\ \  \vect{v}^+=(v_2,v_3,v_4),$$ and define its
\emph{quadratic defect} by
\[
\vect{B}(\vect{v})\eqdef\vect{v}^-\wedge\vect{v}^+\in\bigwedge^2\RR^3.
\]
Let $\rho_2(\tens{H})$ be the rank of its two-against-two flattening.
Equivalently, if $L$ is the degree-twelve moment functional of
$\tens{H}$, then $\rho_2(\tens{H})$ is the rank of the bilinear form $(f,g)\mapsto L(fg)$ on $\set{P}_6$.  
\begin{lemma}[Quadratic defect bound]
\label{lem:quartic-defect-budget}
Let $\tens{H}\in\set{H}_4(4)$ have a real symmetric decomposition
$\tens{H}=\sum_{i=1}^N \gamma_i \vect{v}_i^{\otimes4},$
with $\gamma_i\in\RR\setminus\{0\}$ and $\vect{v}_i\in\RR^4$.  If
$d_{\mathrm{def}}=\dim\Span_i\vect{B}(\vect{v}_i)$, then
\begin{equation*}
\rho_2(\tens{H})+2d_{\mathrm{def}}\le N.
\end{equation*}
\end{lemma}

\begin{proof}
Let $\set{C}\subseteq\RR^N$ be the space obtained by evaluating all
quadratic forms on the points $\vect{v}_i$, and let $\set{E}\subseteq\RR^N$ be the
row space of the bivector columns $\vect{B}(\vect{v}_i)$.  The coordinates of
$\vect{B}(\vect{v})$ are quadratic forms, so $\set{E}\subseteq \set{C}$ and $\dim \set{E}=d_{\mathrm{def}}$.
Equip $\RR^N$ with the nondegenerate diagonal pairing
$\langle\vect{x},\vect{y}\rangle_\gamma\eqdef\sum_{i=1}^N \gamma_ix_iy_i.$ 
For a subspace $\set{C}\subseteq\RR^N$, write
$$\set{C}^{\perp_\gamma}\eqdef
\{\vect{x}\in\RR^N:\langle\vect{x},\vect{c}\rangle_\gamma=0
\text{ for every }\vect{c}\in\set{C}\}.$$
The adjacent Hankel identities applied to the displayed decomposition give
\[
\sum_{i=1}^N \gamma_i \vect{v}_i^{\otimes2}\otimes \vect{B}(\vect{v}_i)=0.
\]
Consequently $\set{E}\subseteq \set{C}\cap \set{C}^{\perp_\gamma}$.  Quadratic evaluation
maps onto $\set{C}$, and the two-against-two flattening is the pullback of
the pairing restricted to $\set{C}$.  Thus its rank is the rank of that
restricted pairing.  Its radical contains $\set{E}$, so
$\rho_2(\tens{H})\le\dim\set{C}-\dim\set{E}$.
Nondegeneracy on the ambient space gives
$$\dim\set{C}^{\perp_\gamma}=N-\dim\set{C}\ge\dim\set{E}.$$
Combining these inequalities with $\dim\set{E}=d_{\mathrm{def}}$ gives
$\rho_2(\tens{H})+2d_{\mathrm{def}}\le N$.
\end{proof}

\begin{proposition}[A symmetric rank lower bound for the central quartic]
\label{prop:quartic-certificates}
The tensor $\tens{H}_*$ defined in \eqref{eq:main-counterexample-tensor}
is concise and satisfies $p_3(\tens{H}_*)\le10$ and $\rho_2(\tens{H}_*)=7$.
Moreover $s(\tens{H}_*)\ge11$.
\end{proposition}

\begin{proof}
Each first-mode coordinate is detected by a triple whose index sum is
$10-i$, so $\tens{H}_*$ is concise.
Its compressed two-against-two flattening is the seven-by-seven
reversal matrix, and \cref{prop:partial-upper} gives $p_3(\tens{H}_*)\le10$.
The moment functional $[t^6]$ on $\set{P}_{12}$ has apolar ideal
$(x^7,y^7)$, which is also the differential apolar ideal of
$F\eqdef x^6y^6$.
Here the \emph{differential apolar ideal} of a binary form $F$ is
$$\{g\in\RR[x,y]:g(\partial_x,\partial_y)F=0\}.$$
By \cite[Proposition~7.7 and the preceding discussion]{blekherman2024almost},
the binary form $F$ has real Waring rank $12$ and almost real rank $10$.
The \emph{almost real rank} is the least degree of an apolar form with at most
one double real root or one nonreal conjugate pair, all remaining
roots being simple and real
\cite[Definition~2.8]{blekherman2024almost}.
A binary form is \emph{split squarefree} if it is a nonzero scalar
times a product of pairwise nonproportional real linear forms.
A split squarefree apolar form must have degree at least $12$
by the real apolarity criterion
\cite[Remark~2.7]{blekherman2024almost}.

Suppose $\tens{H}_*=\sum_{i=1}^N\gamma_i\vect{v}_i^{\otimes4}$
with $N\le10$ and all coefficients and directions nonzero.
By \cref{lem:quartic-defect-budget}, the span of the bivectors
$\vect{B}(\vect{v}_i)$ has dimension at most one.
A zero bivector means that its direction is a real projective node.
Separate these $n_0$ node powers and call the remaining subsum
$\tens{G}$; it is Hankel because each separated node power is Hankel.
Each \emph{linear node factor} is a binary linear form vanishing at the
corresponding projective node.
Let $\ell$ be the least common multiple of these factors,
taking $\ell=1$ if $n_0=0$.
If $\tens{G}=0$, then $\ell$ is split squarefree and apolar to $F$
with $\deg\ell\le n_0\le10$, a contradiction.
Otherwise all remaining bivectors are nonzero and proportional,
so their supports are one plane in $\RR^3$.
Choose a nonzero covector $(q_0,q_1,q_2)$ annihilating this plane
and put $q\eqdef q_0x^2+q_1xy+q_2y^2$.
Every remaining direction satisfies
\[
q_0v_1+q_1v_2+q_2v_3=0;
\ \ q_0v_2+q_1v_3+q_2v_4=0.
\]
Its moment functional $L_G$ therefore annihilates both
$(xq)\set{P}_9$ and $(yq)\set{P}_9$.
Since $x\set{P}_9+y\set{P}_9=\set{P}_{10}$, the form $q$ is apolar to $\tens{G}$.
Choose a common real affine chart preserving the degrees of $q$
and the separated node factors.
By \cref{lem:apolar-quotient}, a divisor $g_0\mid q$ gives a Frobenius
presentation of $\tens{G}$.
The cubic window maps onto this quotient because $\deg g_0\le2$.
Thus \cref{thm:complete-rank} gives
$$\Gamma_4(g_0)=\rk(\tens{G})\le N-n_0.$$

The form $a\eqdef\lcm(g_0,\ell)$ is a multiple of both $g_0$ and $\ell$,
so it is apolar to $\tens{G}$ and to the sum of the separated node powers,
and hence to $\tens{H}_*$.
Return to the original homogeneous coordinates; apolarity,
degrees, and real root multiplicities are preserved.
If $g_0$ is split squarefree, so is $a$, and
$$\deg a\le n_0+\deg g_0=n_0+\Gamma_4(g_0)\le N\le10,$$
contrary to the real rank of $F$.
Otherwise $g_0$ is a double real linear factor or an irreducible
quadratic, so $\Gamma_4(g_0)=4$ and $n_0\le N-4\le6$.
Then $a$ has almost real roots and
$\deg a\le n_0+2\le8$, contrary to the almost real rank of $F$.
Hence $s(\tens{H}_*)\ge11$.
\end{proof}

\subsection{An eleven-power decomposition and exact ranks}
Define four elementary polynomials and a polynomial curve by
\begin{align*}
\phi_1(t)\eqdef t,\ \  \phi_2(t)\eqdef t^2-1,\ \
\phi_3(t)\eqdef t^3-2t,\ \  \phi_4(t)\eqdef t^4-3t^2+1.
\end{align*}
Set $\vect{u}(t)\eqdef(\phi_4(t),\phi_3(t),\phi_2(t),\phi_1(t))$.
Put
\begin{equation}\label{eq:quartic-eleven-node-polynomial}
\psi(t)\eqdef t(t^2-2)
\left(t^4-3t^2+\frac12\right)
\left(t^4-3t^2+\frac32\right).
\end{equation}
Its roots are the eleven distinct real numbers
\[
0,\  \pm\sqrt2,\ 
\pm\sqrt{\frac{3+\sqrt3}{2}},\ 
\pm\sqrt{\frac{3-\sqrt3}{2}},\ 
\pm\sqrt{\frac{3+\sqrt7}{2}},\ 
\pm\sqrt{\frac{3-\sqrt7}{2}}.
\]

\begin{proposition}[An elementary eleven-power decomposition]
\label{prop:quartic-eleven-explicit-powers}
The tensor $\tens{H}_*$ defined in \eqref{eq:main-counterexample-tensor}
satisfies the following exact identity, with $\vect{u}$ and $\psi$ defined above:
\begin{equation}\label{eq:quartic-eleven-explicit-powers}
\tens{H}_* =\sum_{\psi(\theta)=0}
\frac{\vect{u}(\theta)^{\otimes4}}{\psi'(\theta)}.
\end{equation}
In particular $s(\tens{H}_*)\le11$. 
\end{proposition}
\begin{proof}
For polynomials $f$, write
$\Lambda(f)\eqdef\sum_{\psi(\theta)=0}\frac{f(\theta)}{\psi'(\theta)}.$ 
Taking the leading coefficient in the Lagrange interpolation
formula gives $\Lambda(t^j)=0$ for $0\le j\le9$ and
$\Lambda(t^{10})=1$. Also $\Lambda$ vanishes on odd polynomials,
since $\psi$ is odd, $\psi'$ is even, and the nodes occur in
opposite pairs. Expanding \eqref{eq:quartic-eleven-node-polynomial} gives
\[
\psi(t)=t^{11}-8t^9+23t^7-28t^5+\frac{51}{4}t^3-\frac32t.
\]
Using $\psi(\theta)=0$ and the preceding low-degree moments gives
\begin{equation}\label{eq:quartic-eleven-interpolation-moments}
\Lambda(t^{12})=8,\ \ 
\Lambda(t^{14})=41,\ \ 
\Lambda(t^{16})=172.
\end{equation}
Indeed these are respectively $8$, $8\cdot8-23$, and
$8\cdot41-23\cdot8+28$.

Take $a_1,\ldots,a_4\in\{1,2,3,4\}$.  The product $\prod_i \phi_{a_i}$
has degree $\sum_i a_i$ and the same parity.  Grouping the products
of even degree at least ten gives
\[
\begin{array}{c|l}
\text{total degree}&\text{terms of degree at least ten}\\ \hline
10&t^{10}\\
12&t^{12}-8t^{10}\\
14&t^{14}-10t^{12}+39t^{10}\\
16&t^{16}-12t^{14}+58t^{12}-144t^{10}.
\end{array}
\]
Here the second line follows from the coefficient
$-(a_i-1)$ of $t^{a_i-2}$ in $\phi_{a_i}$ (zero when $a_i=1$).
For total degree fourteen the only unordered degree lists are $(4,4,4,2)$
and $(4,4,3,3)$, and both give the displayed coefficients.
For total degree sixteen the product is $\phi_4^4$. Thus all four lines follow
by multiplication of the displayed polynomials of degree at most four.
Together with \eqref{eq:quartic-eleven-interpolation-moments}, this gives
\[
\Lambda\left(\prod_i \phi_{a_i}\right)=
\begin{cases}1,&\sum_i a_i=10,\\0,&\sum_i a_i\ne10.\end{cases}
\]
For completeness, the other three even values give
$$8-8=0,\ \  41-80+39=0,\ \  \textrm{and}\ \  172-492+464-144=0;$$
products of odd degree or degree at most eight vanish already. Substituting $a_j=5-i_j$
proves every entry of \eqref{eq:quartic-eleven-explicit-powers},
since $\sum a_j=10$ exactly when $\sum i_j=10$.
\end{proof}


\begin{proof}[Proof of the single-block part of
\cref{thm:main-counterexample}]
\Cref{prop:quartic-certificates} gives $s(\tens{H}_*)\ge11$.
The identity in \cref{prop:quartic-eleven-explicit-powers}
gives the matching upper bound.
If $\rk(\tens{H}_*)\le9$, then \cref{thm:main-comon} with $M=10$
would give $s(\tens{H}_*)=\rk(\tens{H}_*)\le9$, a contradiction.
Thus 
$$10\le\rk(\tens{H}_*)\le p_2(\tens{H}_*)
\le p_3(\tens{H}_*)\le10,$$ where the last inequality is
\cref{prop:partial-upper}.
\end{proof}

\section{Direct sums and unbounded Comon rank gaps}
\label{sec:quartic-direct-sum-amplification}

Let $\set{W}=\bigoplus_{a=1}^k \set{W}_a$, with $\set{W}_a\simeq\RR^4$, and let
$\tens{H}_*^{(k)}=\bigoplus_{a=1}^k \tens{H}_{*,a}$. 
Here $\tens{H}_{*,a}$ is a copy of \eqref{eq:main-counterexample-tensor}
on $\set{W}_a$; thus $\tens{H}_*^{(k)}$ is the tensor defined in
\eqref{eq:main-counterexample-tensor-k}. 
The symmetric rank lower bound follows from a single-block estimate,
a bound on mixed subdecompositions, and a blockwise count. The
CP and partially symmetric ranks are then determined using
\cref{thm:main-additivity} and the ten-term single-block decomposition.

\subsection{Single-block estimates}

\begin{lemma}[A proper support needs at least four additional powers]
\label{lem:direct-sum-proper-support-four-powers}
Let $\tens{H}_*$ be defined in \eqref{eq:main-counterexample-tensor}.
If $\tens{H}_*=\tens{G}+\sum_{i=1}^p \gamma_i \vect{v}_i^{\otimes4}$,
where $p\in\mathbb Z_{\ge0}$, $\gamma_i\in\RR$, $\vect{v}_i\in\RR^4$, and
$\tens{G}\in\Sym^4(\RR^4)$ is supported on a proper subspace of $\RR^4$,
then $p\ge4$.
\end{lemma}
\begin{proof}
Choose a nonzero first-mode covector $\ell$ annihilating the support
of $\tens{G}$. In the moment coordinates of $\tens{H}_*$, covectors are polynomials
of degree at most three. The cubic obtained by contracting $\ell$
is concise in all four variables. In fact a first-mode kernel
element $f\in \RR[t]_{\le 3}$ would satisfy
\[
[t^6](\ell f q_1 q_2)=0\ \ (q_1,q_2\in \RR[t]_{\le 3}).
\]
The products $q_1q_2$ span $\RR[t]_{\le 6}$, and the pairing
$(f,g)\mapsto[t^6](fg)$ is nondegenerate on $\RR[t]_{\le 6}$.
Consequently $\ell f=0$, and hence $f=0$.
Thus this cubic has first flattening rank four. The displayed
decomposition, after contraction by $\ell$, gives at most $p$
cubic powers, so $p\ge4$.
\end{proof}

\begin{lemma}[A dimension bound for a single block]
\label{lem:direct-sum-single-block-charge}
Let $\tens{H}_*$ be defined in \eqref{eq:main-counterexample-tensor}, and let
$b,p\in\mathbb Z_{\ge0}$. Suppose
\[
\tens{H}_*=\sum_{j=1}^{b} \tens{G}_j+\sum_{i=1}^p \gamma_i \vect{v}_i^{\otimes4},
\]
where $0\ne\tens{G}_j\in\Sym^4(\set{U}_j)$ has intrinsic support
$\set{U}_j\subseteq\RR^4$ of dimension $e_j\in\{2,3,4\}$,
$\gamma_i\in\RR$, and $\vect{v}_i\in\RR^4$. If $b\ge1$, then
$p+\sum_{j=1}^{b}(3e_j-2)\ge10.$
If $b=0$, then $p\ge11$.
\end{lemma}
\begin{proof}
If $b=0$, then $p\ge s(\tens{H}_*)=11$. 
Suppose $b\ge1$.
If every $e_j=2$, each $\tens{G}_j$ is a real binary quartic and
can be replaced by at most $3e_j-2=4$ powers, by
\cite[Lemma~3.11]{zhang2016comon} over $\RR$ with $m=4$.
Since $s(\tens{H}_*)=11$, this gives
$p+\sum_j(3e_j-2)\ge11$.
If some $e_j=4$, then $3e_j-2=10$, which already gives the bound.
It remains to consider a three-dimensional term and no
four-dimensional term.  If $p+\sum_j(3e_j-2)\le9$, then
$b=1$, $e_1=3$, and $p\le2$.  This contradicts
\cref{lem:direct-sum-proper-support-four-powers} applied
to $\tens{G}_1$.
Therefore $p+\sum_j(3e_j-2)\ge10$.
\end{proof}

\subsection{Mixed subdecompositions and a block operator}
The next lemma bounds the lengths of mixed subsums.

\begin{lemma}[The cost of a mixed subdecomposition]
\label{lem:direct-sum-mixed-subsource-cost}
Let $\set{W}=\bigoplus_{a=1}^k\set{W}_a$ with $\set{W}_a\simeq\RR^4$,
let $\tens{F}\in\bigoplus_a\Sym^4(\set{W}_a)$, and start with a
shortest real symmetric decomposition of $\tens{F}$.
A power is called \emph{mixed} if its direction belongs to no single $\set{W}_a$.
Its terms partition
into powers supported in individual $\set{W}_a$ and nonzero mixed subsums
$\tens{F}_1,\ldots,\tens{F}_{b}$ with lengths $N_1,\ldots,N_{b}$. Each mixed subsum
has a direct-sum expression
$\tens{F}_j=\bigoplus_{a\in \set{I}_j}\tens{G}_{ja},$ 
where $0\ne\tens{G}_{ja}\in\Sym^4(\set{U}_{ja})$ has intrinsic support
$\set{U}_{ja}\subseteq \set{W}_a$ of dimension $e_{ja}$.  Moreover
\begin{equation*}
|\set{I}_j|\ge2\ \ \text{and}\ \ 
N_j\ge3\sum_{a\in \set{I}_j}e_{ja}-2.
\end{equation*}
\end{lemma}
\begin{proof}
Choose distinct real numbers $\eta_a$ and let $\matr{D}$ act by
$\eta_a$ on $\set{W}_a$.
For a symmetric fourth-order tensor, equality of the actions
$\matr{D}_1,\ldots,\matr{D}_4$ forces its four block indices to agree.
Thus $\set{K}_{\matr{D}}=\bigoplus_a\Sym^4(\set{W}_a)$.

Apply \cref{lem:operator-balanced-components} to the given shortest
symmetric decomposition with $q=4$ and $\set{E}=\set{W}$.
An eigenvector of $\matr{D}$ lies in a single $\set{W}_a$.
Hence zero columns are precisely the individual block powers,
and every nonzero-column component consists of mixed powers.
Its subsum $\tens{F}_j$ lies in $\set{K}_{\matr{D}}$ and satisfies
$s(\tens{F}_j)=N_j$.
For a block-diagonal tensor, contraction by three covectors
supported in one block isolates that block's first support.
Consequently its intrinsic support is
$\set{U}_j\eqdef\bigoplus_{a\in\set{I}_j}\set{U}_{ja}$.
The component directions span $\set{U}_j$ by
\cref{lem:shortest-subsum-support}; since they are mixed,
$|\set{I}_j|\ge2$.
The bound in \cref{lem:operator-balanced-components} now gives
$$N_j\ge3\dim\set{U}_j-2=3\sum_{a\in\set{I}_j}e_{ja}-2.$$
The proof is complete. 
\end{proof}

\subsection{Exact ranks and Comon rank gaps}
\begin{proposition}[Exact real symmetric ranks of all repeated direct sums]
\label{thm:quartic-direct-sum-symmetric-rank}
For every integer $k\ge1$, let $\tens{H}_*^{(k)}$ be the direct sum of $k$
copies of $\tens{H}_*$ as in \eqref{eq:main-counterexample-tensor-k}. Then
\begin{equation*}
s(\tens{H}_*^{(k)})=11k.
\end{equation*}
\end{proposition}
\begin{proof}
The eleven-power identity in \cref{prop:quartic-eleven-explicit-powers}
gives the upper bound.
Apply \cref{lem:direct-sum-mixed-subsource-cost} to a shortest symmetric
decomposition and let $p_a$ count the individual powers in $\set{W}_a$.
Since $|\set{I}_j|\ge2$, each mixed subsum satisfies
$N_j\ge3\sum_{a\in\set{I}_j}e_{ja}-2
\ge\sum_{a\in\set{I}_j}(3e_{ja}-1).$

For a block incident with a mixed subsum, put
$C_a\eqdef p_a+\sum_{j:a\in\set{I}_j}(3e_{ja}-2)$.
In its restricted identity, regard each term with $e_{ja}=1$
as an individual power; this leaves $C_a$ unchanged.
If a term of support dimension at least two remains,
\cref{lem:direct-sum-single-block-charge} gives $C_a\ge10$.
Otherwise the identity consists of $C_a$ signed powers, so $C_a\ge11$.
In either case, the nonempty incidence set gives
$p_a+\sum_{j:a\in\set{I}_j}(3e_{ja}-1)\ge C_a+1\ge11.$ 
For a block with no mixed subsum, we have $p_a\ge11$.
Summing these blockwise estimates yields
\begin{align*}
s(\tens{H}_*^{(k)})=\sum_a p_a+\sum_j N_j
\ge\sum_a\left(p_a+\sum_{j:a\in\set{I}_j}(3e_{ja}-1)\right)
\ge11k.
\end{align*}
\end{proof}


\begin{proof}[Proof of the direct-sum part of
\cref{thm:main-counterexample}]
By \cref{thm:main-additivity}, the CP rank of $\tens{H}_*$ is additive
under direct sums with every real fourth-order tensor.  Induction therefore gives
$\rk(\tens{H}_*^{(k)})=10k$.  Concatenating the ten-term decompositions defining $p_3$ gives the
matching upper bounds for both partial ranks, and CP rank gives
their lower bounds.  Concision follows from the independent
four-dimensional first supports of the copies.  The symmetric rank
identity is \cref{thm:quartic-direct-sum-symmetric-rank}.  Thus
\eqref{eq:main-direct-sum-gap} holds for every $k\ge1$.
\end{proof}

\section{Conclusions}
\label{sec:outlook}

In this paper, we prove that CP rank is additive for the direct sum of a real Hankel
tensor and an arbitrary real tensor of the same order. We also identify
CP rank regimes in which Comon rank equality holds for real Hankel
tensors. At maximal CP rank, our fourth-order example shows that CP rank
and symmetric rank can differ. For third-order real Hankel tensors of maximal CP rank $M$,
our argument leaves Comon rank equality unresolved in general.
More generally, we ask the following question. 

\begin{problem}\label{prob:maximal-comon-other-orders}
Let $\tens{H}\in\set{H}_m(n)$, where $m\ge3$, $m\ne4$, and $n\ge3$.
Let $M$ be as in \eqref{eq:M_value}.
Suppose that $\rk(\tens{H})=M$.
Does $s(\tens{H})=\rk(\tens{H})$ always hold?
\end{problem}

\medskip
\noindent\textbf{Acknowledgments.} 
The author acknowledges the assistance of ChatGPT 6 Astra throughout the mathematical development and preparation of this manuscript. 
The author takes full responsibility for the content of this paper. 
This work was supported by the Guangdong Basic and Applied Basic
Research Foundation (2026A1515050001, 2026A1515010188).

\bibliographystyle{siamplain}
\bibliography{ReferenceTensor}

\end{document}